# Elementary observations
# on Rogers-Szegö polynomials

Johann Cigler

Fakultät für Mathematik

Universität Wien

johann.cigler@univie.ac.at

**Abstract**

The Rogers-Szegö polynomials are natural $q-$analogues of Newton binomials. In general they have no closed expression. We consider some exceptional cases which are products of a factor with a closed formula and another one with nice values for $q = 1$ and $q = -1$. These are related to "normalized" Rogers-Szegö polynomials which have nice expansions and Hankel determinants. The exposition of the paper is elementary and almost self-contained.

## 0. Introduction

We make some elementary observations on the Rogers-Szegö polynomials
$r_n(s,q) = \sum_{j=0}^{n} \begin{bmatrix} n \\ j \end{bmatrix}_q s^j$. For arbitrary $s$ there is no closed expression for $r_n(s,q)$. Notable exceptions are Gauss's identity $r_{2n}(-1,q) = (q;q^2)_n$ and the formulae $r_n(-q,q) = (q;q^2)_{\left\lfloor \frac{n+1}{2} \right\rfloor}$ and $r_n(q,q^2) = (-q;q)_n$.

B. A. Kupershmidt [9] has obtained some generalizations of these results which were the starting point of my investigation (cf. also [3]).

It turns out that for non-negative integers $m$ the polynomials $r_{2n}(-q^m, q)$, $r_{2n+1}(-q^m, q)$, and $r_n(q^{2m+1}, q^2)$ have the factorizations $r_{2n}(-q^m, q) = (q;q^2)_n T_m(1, -q^{2n}, q)$, $r_{2n-1}(-q^m, q) = (q;q^2)_n U_{m-1}(1, -q^{2n}, q)$, and $r_n(q^{2m+1}, q^2) = (-q;q)_n V_m(1, -q^n, -q)$, where $T_m(x,s,q)$, $U_m(x,s,q)$ and $V_m(x,s,q)$ are $q-$Chebyshev polynomials.

Computer experiments suggest that $r_n(-q^m, q^{2k})$, $r_n(-q^m, q^p)$, and $r_n(q^{2m+1}, q^{2p})$, where $p$ is an odd prime, also have factorizations into a factor with closed expression and another one with nice values for $q = \pm 1$. These results still wait for a proof. The values of the corresponding polynomials for $q=\pm 1$ depend on the limits $\lim_{q\to 1}\frac{r_{2n}(-q^m,q^k)}{(-q;q)_{2n}}=k^n$ and $\lim_{q\to 1}\frac{r_{2n+1}(-q^m,q^k)}{(-q;q)_{2n+1}}=mk^n$ and similar results of this type.

More generally I consider "normalized" Rogers-Szegö polynomials

$$f(n,s,q)=\frac{\sum_{j=0}^{n}s^j\begin{bmatrix}n\\j\end{bmatrix}_{q^2}}{(-q;q)_n}$$

and

$$F(n,s,q)=\frac{\sum_{j=0}^{n}(-s)^j\begin{bmatrix}n\\j\end{bmatrix}_q}{(q;q^2)_{\lfloor\frac{n+1}{2}\rfloor}}.$$

Surprisingly the Hankel determinants $d(n,s,q)=\det\left(f(i+j,s,q)\right)_{i,j=0}^n$ and $D(n,s,q)=\det\left(F(i+j,s,q)\right)_{i,j=0}^n$ which are polynomials in $\mathbb{Q}(q)[s]$ have closed formulae. Some of these were already implicitly contained in [7]. My conjecture for $D(n,s,q)$ has recently been proved by D. Stanton [14]. The roots of $D(n,s,q)$ are of the form $q^{\pm m}$ and the roots of $d(n,s,q)$ are of the form $q^{\pm(2m+1)}$ for integers $m\in\mathbb{Z}$. Thus the sequence $(D(n,s,q))_{n\geq 0}$ vanishes for all large $n$ if and only if $s=q^m$ and the sequence $(d(n,s,q))_{n\geq 0}$ vanishes for all large $n$ if and only if $s=q^{2m+1}$ for some non-negative integer $m$. So these determinants show again the special role of the Rogers-Szegö polynomials $r_n(-q^m,q)$ and $r_n(q^{2m+1},q^2)$.

I want to thank Ole Warnaar, Dennis Stanton and Michael Schlosser for useful comments and the references [1] and [7].

**1. Background material**

**1.1. q-commuting operators**

Let $q$ be a real number with $|q|<1$. As usual we set $[x]=[x]_q=\frac{1-q^x}{1-q}$ for real numbers $x$ and let $[n]!=\prod_{j=1}^{n}[j]$ for $n\in\mathbb{N}.$ The $q-$binomial coefficients $\begin{bmatrix}x\\k\end{bmatrix}$ are defined as

$$\begin{bmatrix} x \\ k \end{bmatrix} = \begin{bmatrix} x \\ k \end{bmatrix}_q = \frac{\prod_{j=0}^{k-1}(1-q^{x-j})}{\prod_{j=0}^{k-1}(1-q^{k-j})}.$$ We shall also use the $q$-Pochhammer symbols

$$(x;q)_n = \prod_{j=0}^{n-1}(1-q^j x) \text{ and } (x;q)_\infty = \prod_{j=0}^{\infty}(1-q^j x).$$

The $q$-binomial coefficients satisfy two recurrence relations

$$\begin{bmatrix} x+1 \\ k \end{bmatrix} = q^k \begin{bmatrix} x \\ k \end{bmatrix} + \begin{bmatrix} x \\ k-1 \end{bmatrix} \tag{1.1}$$

and

$$\begin{bmatrix} x+1 \\ k \end{bmatrix} = \begin{bmatrix} x \\ k \end{bmatrix} + q^{x+1-k} \begin{bmatrix} x \\ k-1 \end{bmatrix}. \tag{1.2}$$

If $A$ and $B$ are linear operators on $\mathbb{C}[x,s]$ which $q$-commute, i.e. satisfy $BA = qAB$, then

$$(A+B)^n = \sum_{k=0}^{n} \begin{bmatrix} n \\ k \end{bmatrix} A^k B^{n-k}. \tag{1.3}$$

This fact has been found by M.P. Schützenberger [15] and has been rediscovered in [2].

The proof is almost trivial since

$$(A+B)(A+B)^n = (A+B)\sum_k \begin{bmatrix} n \\ k \end{bmatrix} A^k B^{n-k} = \sum_k \left(\begin{bmatrix} n \\ k-1 \end{bmatrix} + q^k \begin{bmatrix} n \\ k \end{bmatrix}\right) A^k B^{n+1-k} = \sum_k \begin{bmatrix} n+1 \\ k \end{bmatrix} A^k B^{n+1-k}.$$

By comparing coefficients (1.3) is equivalent with the $q$-exponential law

$$\exp(A+B) = \exp(A)\exp(B) \tag{1.4}$$

for the first $q$-exponential series

$$\exp(z) = \exp_q(z) = \sum_{n \geq 0} \frac{z^n}{[n]!}. \tag{1.5}$$

All considered series will be interpreted as formal power series.

Let me first recall some well-known facts from this point of view (cf. e.g. [2]).

Let $\underline{x}$ and $\underline{s}$ be the multiplication operators with $x$ and $s$ respectively, and let $\varepsilon$ be the augmenting operator defined by $\varepsilon f(x,s) = f(x,qs)$ for $f \in \mathbb{C}[x,s]$.

Then $(A,B) = (\underline{s}\varepsilon, \varepsilon)$, $(A,B) = (\underline{s}, \underline{x}\varepsilon)$ and $(A,B) = (\underline{s}, -\underline{s}\varepsilon)$ are some $q$-commuting operators which will be used in the sequel.

For example $A = \underline{s}$ and $B = \underline{x}\varepsilon$ are $q$-commuting because
$BAf(x,s) = \underline{x}\varepsilon \underline{s} f(x,s) = qsxf(x,qs) = q\underline{s}\underline{x}\varepsilon f(x,s) = qBAf(x,s)$.

Therefore we get $(\underline{s} + \underline{x}\varepsilon)^n = \sum_{k=0}^{n} \begin{bmatrix} n \\ k \end{bmatrix} (\underline{s})^k (\underline{x}\varepsilon)^{n-k}$. If we apply this operator identity to the constant polynomial 1 we get the bivariate Rogers-Szegö polynomials

$$r_n(x,s,q) = \sum_{k=0}^{n} \begin{bmatrix} n \\ k \end{bmatrix}_q s^k x^{n-k} = \sum_{k=0}^{n} \begin{bmatrix} n \\ k \end{bmatrix}_q (\underline{s})^k (\underline{x}\varepsilon)^{n-k} 1 = (\underline{s} + \underline{x}\varepsilon)^n 1. \qquad (1.6)$$

Note that $r_n(1,s,q) = r_n(s,q)$.

There is also a second $q$-exponential series

$$\mathrm{Exp}(z) = \mathrm{Exp}_q(z) = \sum q^{\binom{n}{2}} \frac{z^n}{[n]!}. \qquad (1.7)$$

Note that

$$\exp(z)\mathrm{Exp}(-z) = 1. \qquad (1.8)$$

To show this let $A = \underline{s}$ and $B = -\underline{s}\varepsilon$. Then $\exp(\underline{s})\exp(-\underline{s}\varepsilon) = \exp(\underline{s}(1-\varepsilon))$. Since $\underline{s}(1-\varepsilon)1 = 0$ we get

$$1 = \exp(\underline{s}(1-\varepsilon))1 = \exp(\underline{s})\exp(-\underline{s}\varepsilon)1 = \exp(\underline{s}) \sum_{n \geq 0} (-1)^n q^{\binom{n}{2}} \frac{s^n}{[n]!} = \exp(s)\mathrm{Exp}(-s).$$

Instead of $\exp(z)$ we will also consider the variant

$$e(z) = e_q(z) = \exp_q\left(\frac{z}{1-q}\right) = \sum_{n \geq 0} \frac{z^n}{(q;q)_n} \qquad (1.9)$$

which is the most convenient one in the context of basic hypergeometric series.

From our point of view the $q$-analogue $\dfrac{\exp(z) - \exp(qz)}{(1-q)z} = \exp(z)$ of $\dfrac{d}{dx} e^x = e^x$ implies the well-known results

$$e_q(qz) = (1-z)e_q(z) \qquad (1.10)$$

and by iteration

$$e_q(z) = \frac{1}{(z;q)_\infty}. \qquad (1.11)$$

Therefore
$$e_{q^2}(z^2) = \frac{1}{\prod_{j\geq 0}(1-q^{2j}z^2)} = \frac{1}{\prod_{j\geq 0}(1+q^j z)\prod_{j\geq 0}(1-q^j z)} = e_q(-z)e_q(z). \quad (1.12)$$

The generating function of the Rogers-Szegö polynomials is
$$\sum_{n\geq 0} r_n(x,s,q)\frac{z^n}{[n]!} = \exp(xz)\exp(sz). \quad (1.13)$$

This can be proved by comparing coefficients or by the exponential version of (1.6).

An equivalent version is
$$\sum_{n\geq 0} \frac{r_n(s,q)}{(q;q)_n} z^n = \frac{1}{(z;q)_\infty (sz;q)_\infty} = e_q(z)e_q(sz). \quad (1.14)$$

Note also
$$(s;q)_n = (-\underline{s}\varepsilon + \varepsilon)^n 1 = \sum_{k=0}^n \begin{bmatrix} n \\ k \end{bmatrix}(-\underline{s}\varepsilon)^k \varepsilon^{n-k} 1 = \sum_{k=0}^n (-1)^k q^{\binom{k}{2}} \begin{bmatrix} n \\ k \end{bmatrix} s^k \quad (1.15)$$

and the equivalent $q-$ binomial theorem for $q-$ series
$$\sum_{n\geq 0} \frac{(s;q)_n}{(q;q)_n} z^n = e_q(-s\varepsilon z)e_q(\varepsilon z)1 = \frac{e_q(z)}{e_q(sz)} = \frac{(sz;q)_\infty}{(z;q)_\infty}. \quad (1.16)$$

We will also need the following $q-$ analogue of $e^{(x+y)^2} = e^{x^2} e^{2xy} e^{y^2}$.

**Lemma 1.1 ([8], Corollary 3.3)**

Let $A$ and $B$ be linear operators on $\mathbb{C}[x,s]$ which satisfy $BA = qAB$.

Then $\exp_{q^2}((A+B)^2) = \exp_{q^2}(A^2)\exp_q((1+q)AB)\exp_{q^2}(B^2)$ or equivalently
$$e_{q^2}\left((A+B)^2\right) = e_{q^2}(A^2)e_q(AB)e_{q^2}(B^2). \quad (1.17)$$

**Proof**

Since $B^n A = Aq^n B^n$ we get $e_q(B)A = Ae_q(qB)$ and therefore
$e_q(B)Ae_q(B)^{-1} = Ae_q(qB)e_q(B)^{-1} = A(1-B)e_q(B)e_q(B)^{-1} = A(1-B).$

Thus $e_q(B)A^n e_q(B)^{-1} = (A(1-B))^n$ and $e_q(B)e_q(A)e_q(B)^{-1} = e_q(A(1-B))$.

This implies

$$e_q(B)e_q(A) = e_q(A - AB)e_q(B) = e_q(A)e_q(-AB)e_q(B). \tag{1.18}$$

Therefore by (1.12)

$$e_{q^2}\left((A+B)^2\right) = e_q(-(A+B))e_q((A+B)) = e_q(-A)e_q(-B)e_q(A)e_q(B)$$
$$= e_q(-A)e_q(A)e_q(AB)e_q(-B)e_q(B) = e_{q^2}(A^2)e_q(AB)e_{q^2}(B^2).$$

### 1.2. The Rogers-Szegö polynomials

The Rogers-Szegö polynomials $r_n(x, s, q)$ satisfy

$$r_n(x, s, q) = (\underline{s} + \underline{x}\varepsilon)r_{n-1}(x, s, q) = sr_{n-1}(x, s, q) + xr_{n-1}(x, qs, q).$$

This implies the well-known recurrence

$$r_n(s, q) = (1 + s)r_{n-1}(s, q) + (q^{n-1} - 1)sr_{n-2}(s, q). \tag{1.19}$$

For

$$r_n(s, q) = (1 + s)r_{n-1}(s, q) + (r_{n-1}(qs, q) - r_{n-1}(s, q))$$
$$= (1 + s)r_{n-1}(s, q) + (q^{n-1} - 1)sr_{n-2}(s, q).$$

The last identity follows from

$$r_{n-1}(qs, q) - r_{n-1}(s, q) = \sum_k \begin{bmatrix} n-1 \\ k \end{bmatrix}(q^k - 1)s^k = (q^{n-1} - 1)\sum_k \begin{bmatrix} n-2 \\ k-1 \end{bmatrix}s^k = sr_{n-2}(s, q)$$

by noting that $\begin{bmatrix} n \\ k \end{bmatrix} = \frac{1 - q^n}{1 - q^k}\begin{bmatrix} n-1 \\ k-1 \end{bmatrix}$.

For later use let us mention the following consequence:

$$r_n(s, q) = \left(1 + (1 + q)q^{n-2}s + s^2\right)r_{n-2}(s, q) - \left(1 - q^{n-3}\right)\left(1 - q^{n-2}\right)s^2 r_{n-4}(s, q). \tag{1.20}$$

**Proof**

By (1.19) we have

$$r_n(s, q) = (1 + s)r_{n-1}(s, q) + (q^{n-1} - 1)sr_{n-2}(s, q),$$
$$r_{n-1}(s, q) = (1 + s)r_{n-2}(s, q) + (q^{n-2} - 1)sr_{n-3}(s, q),$$
$$(1 + s)r_{n-3}(s, q) = r_{n-2}(s, q) + (1 - q^{n-3})sr_{n-4}(s, q).$$

This implies

$$r_n(s,q) = (1+s)\big((1+s)r_{n-2}(s,q) + (q^{n-2}-1)sr_{n-3}(s,q)\big) + (q^{n-1}-1)sr_{n-2}(s,q)$$
$$= \big((1+s)^2 + (q^{n-1}-1)s\big)r_{n-2}(s,q) + (q^{n-2}-1)s\big(r_{n-2}(s,q) + (1-q^{n-3})sr_{n-4}(s,q)\big)$$
$$= \big(1 + (1+q)q^{n-2}s + s^2\big)r_{n-2}(s,q) - (1-q^{n-3})(1-q^{n-2})s^2 r_{n-4}(s,q).$$

For $x=1$ and $s=-1$ (1.19) reduces to

$r_n(-1,q) = (1-q^{n-1})r_{n-2}(-1,q)$. Since $r_0(-1,q)=1$ and $r_1(-1,q)=0$

and thus to

**Gauss's identity**

$$\sum_{k=0}^{2n+1}(-1)^k \begin{bmatrix}2n+1\\k\end{bmatrix} = 0,$$
$$\sum_{k=0}^{2n}(-1)^k \begin{bmatrix}2n\\k\end{bmatrix} = (q;q^2)_n. \tag{1.21}$$

This identity is equivalent with

$$\exp(z)\exp(-z) = \sum_{n\geq 0}(q;q^2)_n \frac{z^{2n}}{[2n]!} \tag{1.22}$$

or equivalently with (1.12).

For $x=1$ and $s=-q$ we get

$$r_n(-q,q) = (q;q^2)_{\left\lfloor \frac{n+1}{2} \right\rfloor}. \tag{1.23}$$

Let us also note the following $q$–analogue of $e^{\frac{x}{2}}e^{\frac{x}{2}} = e^x$:

$$\exp_{q^2}\left(\frac{z}{[2]}\right)\exp_{q^2}\left(\frac{qz}{[2]}\right) = \exp_q(z). \tag{1.24}$$

or equivalently

$$e_{q^2}(qz)e_{q^2}(z) = e_q(z) \tag{1.25}$$

which is obvious by (1.11).

By comparing coefficients this is equivalent with

$$r_n(q,q^2) = \sum_{k=0}^n q^k \begin{bmatrix}n\\k\end{bmatrix}_{q^2} = (1+q)(1+q^2)\cdots(1+q^n) = (-q;q)_n. \tag{1.26}$$

A. Berkovich and S.O. Warnaar [1] have found generating functions for $r_{2n}(s,q)$ and $r_{2n+1}(s,q)$ and related facts. I will recall these results with a different proof.

**Lemma 1.2 ([1], (8.9) and (8.11))**

$$\sum_{n \geq 0} \frac{r_{2n}(s,q)}{(q^2;q^2)_n} z^n = \frac{1}{e_q(-sz)} e_{q^2}(s^2 z) e_{q^2}(z) \tag{1.27}$$

$$\sum_{n \geq 0} \frac{r_{2n+1}(s,q)}{(q^2;q^2)_n} z^n = (1+s) \frac{1}{e_q(-qsz)} e_{q^2}(s^2 z) e_{q^2}(z). \tag{1.28}$$

**Proof**

By (1.17) and (1.13) we get

$$\sum_{n \geq 0} \frac{r_{2n}(s,q)}{(q^2;q^2)_n} z^n = e_{q^2}\left(( \underline{s}+\varepsilon)^2 z \right)1 = e_{q^2}(\underline{s}^2 z) e_q(\underline{s}\varepsilon z) e_{q^2}(\varepsilon^2 z)1 = e_{q^2}(s^2 z) \frac{1}{e_q(-sz)} e_{q^2}(z).$$

For the second identity observe that by (1.10)

$$\sum_{n \geq 0} \frac{r_{2n+1}(s,q)}{(q^2;q^2)_n} z^n = (\underline{s}+\varepsilon) e_{q^2}\left((\underline{s}+\varepsilon)^2 z\right)1 = (\underline{s}+\varepsilon) e_{q^2}(\underline{s}^2 z) e_q(\underline{s}\varepsilon z) e_{q^2}(\varepsilon^2 z)1$$

$$= (\underline{s}+\varepsilon) e_{q^2}(s^2 z) \frac{1}{e_q(-sz)} e_{q^2}(z) = s e_{q^2}(s^2 z) \frac{1}{e_q(-sz)} e_{q^2}(z) + e_{q^2}(q^2 s^2 z) \frac{1}{e_q(-qsz)} e_{q^2}(z)$$

$$= s(1+sz) e_{q^2}(s^2 z) \frac{1}{e_q(-qsz)} e_{q^2}(z) + \left(1 - s^2 z\right) e_{q^2}(s^2 z) \frac{1}{e_q(-qsz)} e_{q^2}(z) = (1+s) \frac{e_{q^2}(s^2 z) e_{q^2}(z)}{e_q(-qsz)}.$$

Comparing coefficients and observing (1.13) we get

**Corollary 1.1 ([1], (8.8))**

$$r_{2n}(s,q) = \sum_{k=0}^{n} (-q;q)_k \, q^{\binom{k}{2}} s^k \begin{bmatrix} n \\ k \end{bmatrix}_{q^2} r_{n-k}(s^2, q^2) \tag{1.29}$$

*and*

$$r_{2n+1}(s,q) = (1+s) \sum_{k=0}^{n} (-q;q)_k \, q^{\binom{k+1}{2}} s^k \begin{bmatrix} n \\ k \end{bmatrix}_{q^2} r_{n-k}(s^2, q^2). \tag{1.30}$$

Another consequence which is a $q-$ analogue of $(1+s)^{2n} = \sum_{k=0}^{n} \binom{n}{k} s^{2k} \left(1 + \frac{1}{s}\right)^k (1+s)^{n-k}$ is

**Corollary 1.2 (A. Berkovich and S.O. Warnaar [1], Theorem 8.1)**

$$r_{2n}(s,q) = \sum_{k=0}^{n} s^{2k} \left(-\frac{q}{s};q^2\right)_k (-s,q^2)_{n-k} \begin{bmatrix} n \\ k \end{bmatrix}_{q^2}, \qquad (1.31)$$

*and*

$$r_{2n+1}(s,q) = \sum_{k=0}^{n} s^{2k} \left(-\frac{q}{s},q^2\right)_k (-s,q^2)_{n+1-k} \begin{bmatrix} n \\ k \end{bmatrix}_{q^2}$$

$$= (1+s)\sum_{k=0}^{n} s^{2k} \left(-\frac{q}{s},q^2\right)_k (-q^2 s,q^2)_{n-k} \begin{bmatrix} n \\ k \end{bmatrix}_{q^2}. \qquad (1.32)$$

**Proof**

By comparing coefficients (1.31) is equivalent with

$$\sum_{n\geq 0} \frac{r_{2n}(s,q)}{(q^2;q^2)_n} z^n = \sum_{n} s^{2n} \left(-\frac{q}{s};q^2\right)_n \frac{z^n}{(q^2;q^2)_n} \sum_{n} (-s;q^2)_n \frac{z^n}{(q^2;q^2)_n}$$

$$= \frac{e_{q^2}(s^2 z)}{e_{q^2}(-qsz)} \frac{e_{q^2}(z)}{e_{q^2}(-sz)} = \frac{e_{q^2}(s^2 z)e_{q^2}(z)}{e_q(-sz)}$$

and (1.32) with

$$\sum_{n\geq 0} \frac{r_{2n+1}(s,q)}{(q^2;q^2)_n} z^n = (1+s)\sum_{n} s^{2n} \left(-\frac{q}{s};q^2\right)_n \frac{z^n}{(q^2;q^2)_n} \sum_{n} (-q^2 s;q^2)_n \frac{z^n}{(q^2;q^2)_n}$$

$$= (1+s)\frac{e_{q^2}(s^2 z)}{e_{q^2}(-qsz)} \frac{e_{q^2}(z)}{e_{q^2}(-q^2 sz)} = (1+s)\frac{e_{q^2}(s^2 z)e_{q^2}(z)}{e_q(-qsz)}.$$

There are curiously similar formulae if we replace $r_{2n}(s,q)$ with $r_n(s,q^2)$.

**Lemma 1.3**

*The generating function of $r_n(s,q^2)$ is given by*

$$\sum_{n} \frac{r_n(s,q^2)}{(q;q)_n} z^n = \frac{e_q(sz)e_q(z)}{e_{q^2}(qsz^2)} = \frac{e_q(sz)}{e_q(-z)} \frac{e_{q^2}(z^2)}{e_{q^2}(qsz^2)}. \qquad (1.33)$$

*Equivalently we get*

$$r_n(s,q^2) = \sum_{j=0}^{\lfloor \frac{n}{2} \rfloor} (-1)^j q^{j^2} (q;q^2)_j \begin{bmatrix} n \\ 2j \end{bmatrix} s^j r_{n-2j}(s,q). \qquad (1.34)$$

*This is again equivalent with*

$$r_n(s,q^2) = \sum_{j=0}^{\left\lfloor \frac{n}{2} \right\rfloor} \begin{bmatrix} n \\ 2j \end{bmatrix} (q;q^2)_j (qs;q^2)_j s^{n-2j} \left(-\frac{1}{s};q\right)_{n-2j} \quad (1.35)$$

*and with*

$$\sum_{k=0}^{2n} (-1)^k (-s;q)_k \begin{bmatrix} 2n \\ k \end{bmatrix} r_{2n-k}(s,q^2) = (q;q^2)_n (qs;q^2)_n,$$

$$\sum_{k=0}^{2n+1} (-1)^k (-s;q)_k \begin{bmatrix} 2n+1 \\ k \end{bmatrix} r_{2n+1-k}(s,q^2) = 0. \quad (1.36)$$

**Proof**

The equivalences are easily verified. By (1.33) and (1.14) we get (1.34).

The right-hand side of (1.33) and (1.16) give (1.35) and (1.36) is equivalent with

$$\frac{e_{q^2}(z^2)}{e_{q^2}(qsz^2)} = \frac{e_q(-z)}{e_q(sz)} \sum_n \frac{r_n(s,q^2)}{(q;q)_n} z^n.$$

Therefore it suffices to show that the right-hand side of (1.35) satisfies the same recursion and the same initial values as the left-hand side. The initial values coincide for $n=0$ and $n=1$ and the recursion can be seen by applying the $q$ – Zeilberger algorithm to the right-hand side. The Mathematica implementation qZeil by P. Paule and A. Riese [10] gives

```
qZeil[qBinomial[n, 2 j, q] qPochhammer[q, q^2, j] qPochhammer[q s, q^2, j] s^(n - 2 j)
  qPochhammer[-1 / s, q, n - 2 j], {j, 0, n}, n, 2]
SUM[n] == - (1 - q^(-2+2 n)) s SUM[-2 + n] + (1 + s) SUM[-1 + n]
```

## 2. Normalized Rogers-Szegö polynomials

### 2.1. A useful expansion

Let us introduce the "normalized" Rogers-Szegö polynomials

$$f(n,s,q) = \frac{\sum_{j=0}^{n} s^j \begin{bmatrix} n \\ j \end{bmatrix}_{q^2}}{(-q;q)_n} \quad (2.1)$$

and

$$F(n,s,q) = \frac{\sum_{j=0}^{n}(-s)^j \begin{bmatrix} n \\ j \end{bmatrix}_q}{(q;q^2)_{\left\lfloor \frac{n+1}{2} \right\rfloor}}. \tag{2.2}$$

which satisfy $f(n,q,q) = F(n,q,q) = 1$.

**Theorem 2.1**

*The normalized Rogers-Szegö polynomials have the following expansions:*

$$f(n,s,q) = \frac{\sum_{j=0}^{n} s^j \begin{bmatrix} n \\ j \end{bmatrix}_{q^2}}{(-q;q)_n} = \sum_{j=0}^{n}(-1)^j \begin{bmatrix} n \\ j \end{bmatrix}_q \frac{\prod_{i=0}^{j-1}(q^{2i+1}-s)}{(-q;q)_j}, \tag{2.3}$$

$$F(2n,s,q) = \frac{\sum_{j=0}^{2n}(-s)^j \begin{bmatrix} 2n \\ j \end{bmatrix}_q}{(q;q^2)_n} = \sum_{k=0}^{n} \begin{bmatrix} n \\ k \end{bmatrix}_{q^2} \frac{\prod_{j=0}^{2k-1}(q^j-s)}{(q;q^2)_k} \tag{2.4}$$

and

$$F(2n+1,s,q) = \frac{\sum_{j=0}^{2n+1}(-s)^j \begin{bmatrix} 2n+1 \\ j \end{bmatrix}_q}{(q;q^2)_{n+1}} = \sum_{k=0}^{n} \begin{bmatrix} n \\ k \end{bmatrix}_{q^2} \frac{\prod_{j=0}^{2k}(q^j-s)}{(q;q^2)_{k+1}}. \tag{2.5}$$

**Proof**

For (2.3) we have to show that

$$f(n,s,q) = \frac{\sum_{j=0}^{n} s^j \begin{bmatrix} n \\ j \end{bmatrix}_{q^2}}{(-q;q)_n} = \sum_{j=0}^{n}(-1)^j q^{j^2} \left(q^{n+1-j};q\right)_j \frac{\left(\frac{s}{q^{2j-1}};q^2\right)_j}{(q^2;q^2)_j} = \sum_{j=0}^{n}(-1)^j \begin{bmatrix} n \\ j \end{bmatrix}_q \frac{\prod_{i=0}^{j-1}(q^{2i+1}-s)}{(-q;q)_j}.$$

To prove this observe that

$$\sum_{n\geq 0}\sum_{j=0}^{n} s^j \begin{bmatrix} n \\ j \end{bmatrix}_{q^2} \frac{z^n}{(q^2;q^2)_n} = e_{q^2}(sz) e_{q^2}(z) = \frac{1}{(sz;q^2)_\infty (z;q^2)_\infty} = \frac{(qz;q^2)_\infty}{(sz;q^2)_\infty} \frac{1}{(z;q^2)_\infty (qz;q^2)_\infty}$$

$$= \sum_k \frac{\left(\frac{q}{s};q^2\right)_k}{(q^2;q^2)_k} (sz)^k \sum_\ell \frac{z^\ell}{(q;q)_\ell} = \sum_k (-1)^k q^{k^2} \frac{\left(\frac{s}{q^{2k-1}};q^2\right)_k}{(q^2;q^2)_k} z^k \sum_\ell \frac{z^\ell}{(q;q)_\ell}$$

implies

$$\sum_{j=0}^{n} s^j \begin{bmatrix} n \\ j \end{bmatrix}_{q^2} = (q^2;q^2)_n \sum_{j=0}^{n} (-1)^j q^{j^2} \left(\frac{s}{q^{2j-1}};q^2\right)_j \frac{1}{(q^2;q^2)_j (q;q)_{n-j}}$$

$$= (-q;q)_n \sum_{j=0}^{n} (-1)^j q^{j^2} \left(\frac{s}{q^{2j-1}};q^2\right)_j \frac{(q;q)_n}{(q^2;q^2)_j (q;q)_{n-j}}.$$

To prove the other identities we compare coefficients in

$$\sum_n \sum_{j=0}^{n} (-s)^j \begin{bmatrix} n \\ j \end{bmatrix} \frac{z^n}{(q;q)_n} = e_q(-sz) e_q(z) = \frac{e_q(-sz)}{e_q(-z)} e_q(-z) e_q(z)$$

$$= \sum_{k\geq 0} q^{\binom{k}{2}} \left(\frac{s}{q^{k-1}};q\right)_k \frac{z^k}{(q;q)_k} \cdot \sum_{\ell \geq 0} \frac{z^{2\ell}}{(q^2;q^2)_\ell}$$

gives

$$\sum_{j=0}^{n} (-s)^j \begin{bmatrix} n \\ j \end{bmatrix} = \sum_{j+2\ell=n} q^{\binom{j}{2}} \left(\frac{s}{q^{j-1}};q\right)_j \frac{(q;q)_n}{(q;q)_j (q^2;q^2)_\ell}.$$

This implies

$$F(2n,s,q) = \sum_{\ell=0}^{n} q^{\binom{2n-2\ell}{2}} \left(\frac{s}{q^{2n-2j-1}};q\right)_{2n-2j} \frac{(q^2;q^2)_n}{(q;q)_{2n-2j} (q^2;q^2)_j},$$

$$F(2n+1,s,q) = \sum_{\ell=0}^{n} q^{\binom{2n-2\ell+1}{2}} \left(\frac{s}{q^{2n-2j}};q\right)_{2n-2j+1} \frac{(q^2;q^2)_n}{(q;q)_{2n-2j+1} (q^2;q^2)_j}.$$

This can be simplified to give (2.4) and (2.5).

## 2.2. Factorizations

### Corollary 2.1 (B. A. Kupershmidt [9])

*For $m \in \mathbb{N}$ the polynomials $r_n(q^{2m+1}, q^2)$ are divisible by $(-q;q)_n$ and the polynomials $r_n(-q^m, q)$ are divisible by $(q;q^2)_{\lfloor \frac{n+1}{2} \rfloor}$. More precisely*

$$f\left(n, q^{2m+1}, q\right) = \frac{\sum_{j=0}^{n} q^{(2m+1)j} \begin{bmatrix} n \\ j \end{bmatrix}_{q^2}}{(-q;q)_n} = \sum_{j=0}^{m} (-1)^j q^{j^2} \begin{bmatrix} m \\ j \end{bmatrix}_{q^2} \left(q^{n-j+1}; q\right)_j \in \mathbb{Z}[q], \qquad (2.6)$$

$$F(2n, q^m, q) = \frac{\sum_{j=0}^{2n} (-1)^j q^{mj} \begin{bmatrix} 2n \\ j \end{bmatrix}_q}{(q;q^2)_n} = \sum_{k=0}^{\lfloor \frac{m}{2} \rfloor} q^{\binom{2k}{2}} \begin{bmatrix} m \\ 2k \end{bmatrix} \left(q^{2n-2k+2}; q^2\right)_k \in \mathbb{Z}[q] \qquad (2.7)$$

*and*

$$F(2n+1, q^m, q) = \sum_{k=0}^{\lfloor \frac{m-1}{2} \rfloor} q^{\binom{2k+1}{2}} \begin{bmatrix} m \\ 2k+1 \end{bmatrix} \left(q^{2n-2k+2}; q^2\right)_k \in \mathbb{Z}[q]. \qquad (2.8)$$

**Proof**

Observe that

$$\sum_{j=0}^{n} (-1)^j \begin{bmatrix} n \\ j \end{bmatrix}_q \frac{\prod_{i=0}^{j-1}\left(q^{2i+1} - q^{2m+1}\right)}{(q;q)_j} = \sum_{j=0}^{n} (-1)^j \frac{(q;q)_n}{(q;q)_j (q;q)_{n-j}} q^{j^2} \frac{\prod_{i=0}^{j-1}\left(1 - q^{2m-2i}\right)}{(q;q)_j}$$
$$= \sum_{j=0}^{m} (-1)^j q^{j^2} \begin{bmatrix} m \\ j \end{bmatrix}_{q^2} \left(q^{n-j+1}; q\right)_j$$

and

$$\sum_{k=0}^{n} \begin{bmatrix} n \\ k \end{bmatrix}_{q^2} \frac{\prod_{j=0}^{2k-1}\left(q^j - q^m\right)}{(q;q^2)_k} = \sum_{k=0}^{n} \frac{(q^2;q^2)_n}{(q^2;q^2)_k (q^2;q^2)_{n-k}} q^{\binom{2k}{2}} \frac{(q;q)_m}{(q;q^2)_k (q;q)_{m-2k}}$$
$$= \sum_{k=0}^{\lfloor \frac{m}{2} \rfloor} q^{\binom{2k}{2}} \begin{bmatrix} m \\ 2k \end{bmatrix} \frac{(q;q)_{2k}}{(q;q^2)_k (q^2;q^2)_k} \left(q^{2n-2k+2}; q^2\right)_k = \sum_{k=0}^{\lfloor \frac{m}{2} \rfloor} q^{\binom{2k}{2}} \begin{bmatrix} m \\ 2k \end{bmatrix} \left(q^{2n-2k+2}; q^2\right)_k.$$

There is also a representation of (2.7) and (2.8) in terms of $q$ − Chebyshev polynomials which sheds new light on some results of [3]. The $q$ − Chebyshev polynomials of the second kind (cf. [5])

$$U_n(x,s,q) = \sum_{k=0}^{\lfloor \frac{n}{2} \rfloor} q^{k^2} \begin{bmatrix} n-k \\ k \end{bmatrix} (-q^{k+1};q)_{n-2k} s^k x^{n-2k} \tag{2.9}$$

satisfy

$$U_m(x,s,q) = (1+q^m) x U_{m-1}(x,s,q) + q^{m-1} s U_{m-2}(x,s,q) \tag{2.10}$$

with initial values $U_{-1}(x,s,q) = 0$ and $U_0(x,s,q) = 1$.

Note that $U_n(x,-1,1) = U_n(x)$ are the classical Chebyshev polynomials of the second kind (cf. OEIS [11], A133156).

The $q$ − Chebyshev polynomials of the first kind (cf. [5]

$$T_n(x,s,q) = \sum_{k=0}^{\lfloor \frac{n}{2} \rfloor} q^{k^2} \frac{[n]}{[n-k]} \begin{bmatrix} n-k \\ k \end{bmatrix} \frac{(-q;q)_{n-1}}{(-q;q)_k (-q^{n-k};q)_k} s^k x^{n-2k} \tag{2.11}$$

satisfy

$$T_m(x,s,q) = (1+q^{m-1}) x T_{m-1}(x,s,q) + q^{m-1} s T_{m-2}(x,s,q) \tag{2.12}$$

with initial values $T_0(x,s,q) = 1$ and $T_1(x,s,q) = x$.

Note that $T_n(x,-1,1) = T_n(x)$ are the classical Chebyshev polynomials of the first kind (cf. OEIS [11], A028297).

Let us define $q$ − Chebyshev polynomials of the third kind $V_n(x,s,q)$ by the recurrence

$$V_m(x,s,q) = (1+q^{2m-1}) x V_{m-1}(x,s,q) - q^{2m-1} s^2 V_{m-2}(x,s,q) \tag{2.13}$$

with initial values $V_0(x,s,q) = 1$ and $V_1(x,s,q) = (1+q)x + qs$.

Note that $V_n(x,-1,1) = V_n(x)$ are the classical Chebyshev polynomials of the third kind (cf. OEIS [11], A228565).

Observe that

$$r_n(q^2 x, q) = (1-qx) r_n(qx,q) + q^{n+1} x r_n(x,q). \tag{2.14}$$

By comparing coefficients this is equivalent with $q^{2k} \begin{bmatrix} n \\ k \end{bmatrix} - q^k \begin{bmatrix} n \\ k \end{bmatrix} = q^{n+1} \begin{bmatrix} n \\ k-1 \end{bmatrix} - q^k \begin{bmatrix} n \\ k-1 \end{bmatrix}$ or $\begin{bmatrix} n \\ k \end{bmatrix} = \frac{1-q^{n+1-k}}{1-q^k} \begin{bmatrix} n \\ k-1 \end{bmatrix} = \frac{[n+1-k]}{[k]} \begin{bmatrix} n \\ k-1 \end{bmatrix}$ which is obviously true.

**Corollary 2.2**

*For $m,n \in \mathbb{N}$ we get*

$$r_{2n}\left(-q^m, q\right) = \left(q;q^2\right)_n T_m\left(1, -q^{2n}, q\right), \qquad (2.15)$$

$$r_{2n-1}\left(-q^m, q\right) = \left(q;q^2\right)_n U_{m-1}\left(1, -q^{2n}, q\right), \qquad (2.16)$$

*and*

$$r_n\left(q^{2m+1}, q^2\right) = \left(-q;q\right)_n V_m\left(1, -q^n, -q\right). \qquad (2.17)$$

To show this let $n$ be fixed and observe that

$$F\left(2n, q^m, q\right) = \left(1 + q^{m-1}\right) F\left(2n, q^{m-1}, q\right) - q^{m-1} q^{2n} F\left(2n, q^{m-2}, q\right)$$

with initial values

$$F\left(2n, q^0, q\right) = \frac{r_{2n}(-1, q)}{\left(q;q^2\right)_n} = 1 \text{ and } F\left(2n, q, q\right) = \frac{r_{2n}(-q, q)}{\left(q;q^2\right)_n} = 1$$

which proves (2.15).

For the second identity we have

$$F\left(2n-1, q^m, q\right) = \left(1 + q^{m-1}\right) F\left(2n-1, q^{m-1}, q\right) - q^{m-1} q^{2n-1} F\left(2n-1, q^{m-2}, q\right)$$

with initial values

$$F\left(2n-1, q^0, q\right) = \frac{r_{2n-1}(-1, q)}{\left(q;q^2\right)_n} = 0 \text{ and } F\left(2n-1, q, q\right) = \frac{r_{2n-1}(-q, q)}{\left(q;q^2\right)_n} = 1$$

which implies $F\left(2n-1, q^m, q\right) = U_{m-1}\left(1, -q^{2n}, q\right)$.

In the same way we get that

$$f\left(n, q^{2m+1}, q\right) = \frac{\sum_{j=0}^{n} q^{(2m+1)j} \begin{bmatrix} n \\ j \end{bmatrix}_{q^2}}{(-q;q)_n} \text{ satisfies}$$

$$f\left(n, q^{2m+1}, q\right) = \left(1 - q^{2m-1}\right) f\left(n, q^{2m-1}, q\right) + q^{2n+2m-1} f\left(n, q^{2m-3}, q\right)$$

with initial values $f\left(n, -\frac{1}{q}, q\right) = \dfrac{\sum_{j=0}^{n} q^{-j} \begin{bmatrix} n \\ j \end{bmatrix}_{q^2}}{(-q;q)_n} = \dfrac{1}{q^n}$ and $f(n, q, q) = 1$.

This implies $f(n, q^3, q) = 1 - q + q^{n+1}$ and therfore $f\left(n, q^{2m+1}, q\right) = V_m\left(1, -q^n, -q\right)$.

To obtain the values for $q \to \pm 1$ note that $T_m(1,-q^{2n},q)_{q=\pm 1} = 1$, $U_m(1,-q^{2n},q)_{q=1} = 1$,

$U_m(1,-q^{2n},q)_{q=-1} = (-1)^{\binom{m+1}{2}} 3^{\lfloor \frac{m}{2} \rfloor}$, $V_m(1,-q^n,-q)_{q=1} = 1$, $V_m(1,-q^{2n},-q)_{q=-1} = 1$,

$V_m(1,-q^{2n+1},-q)_{q=1} = 2m+1$.

**Remarks**

Formulae (2.6) - (2.8) also give explicit formulae for negative $m$ because the change $j \to n-j$ does not change the sums. Thus the alternating sums $r_n(-q^m, q)$ can be explicitly computed for all integers $m \in \mathbb{Z}$ and the non-alternating sums $r_n(q^r, q)$ can be explicitly computed for all half-integers $r \in \frac{1}{2} + \mathbb{Z}$.

For $n \to \infty$ (2.6) converges to $\dfrac{1}{(-q;q)_\infty} \sum_{j \geq 0} \dfrac{q^{(2m+1)j}}{(q^2;q^2)_j} = \sum_{j=0}^{m} (-1)^j q^{j^2} \begin{bmatrix} m \\ j \end{bmatrix} = (q;q^2)_m$.

Formula (2.7) converges to $\dfrac{1}{(q;q^2)_\infty} \sum_{j \geq 0} (-1)^j \dfrac{q^{mj}}{(q;q)_j} = \sum_{k=0}^{\lfloor \frac{m}{2} \rfloor} q^{\binom{2k}{2}} \begin{bmatrix} m \\ 2k \end{bmatrix} = (-q;q)_{m-1}$

and so does (2.8).

**Corollary 2.3**

Let $m \in \mathbb{N}$. Then $\sum_{j=0}^{n} (-1)^j q^{mj} \begin{bmatrix} n \\ j \end{bmatrix}_{q^2}$ is divisible by $(q;q^2)_{\lfloor \frac{n+1}{2} \rfloor}$.

If we choose in (2.3) $s = -q$ then we get $\sum_{j=0}^{n} (-1)^j q^{(2m+1)j} \begin{bmatrix} n \\ j \end{bmatrix}_{q^2}$ is divisible by $(q;-q)_n$.

Since $(q;-q)_n = (q;q^2)_{\lfloor \frac{n+1}{2} \rfloor} (-q^2;q^2)_{\lfloor \frac{n}{2} \rfloor}$ we see that $\sum_{j=0}^{n} (-1)^j q^{(2m+1)j} \begin{bmatrix} n \\ j \end{bmatrix}_{q^2}$ is also divisible by $(q;q^2)_{\lfloor \frac{n+1}{2} \rfloor}$.

On the other hand we know that $\sum_{j=0}^{n} (-1)^j q^{2mj} \begin{bmatrix} n \\ j \end{bmatrix}_{q^2}$ is divisible by

$(q^2;q^4)_{\lfloor \frac{n+1}{2} \rfloor} = (q;q^2)_{\lfloor \frac{n+1}{2} \rfloor} (-q;q^2)_{\lfloor \frac{n+1}{2} \rfloor}$.

Computer experiments suggest more generally the following conjectures.

**Conjecture 2.1**

Let $m, k \in \mathbb{N}$. Then

$$r_n\left(-q^m, q^{2^k}\right) = \sum_{j=0}^{n}(-1)^j q^{mj} \begin{bmatrix} n \\ j \end{bmatrix}_{q^{2^k}} \text{ is divisible by } \left(q; q^2\right)_{\left\lfloor \frac{n+1}{2} \right\rfloor} \text{ in } \mathbb{Z}[q].$$

**Conjecture 2.2**

For a prime number $p$ let $v_p(x)$ denote the $p$-adic valuation of $x$, which is the largest non-negative integer $m$ such that $p^m$ divides $x$ and let $V_p(x) = p^{v_p(x)}$. For integers $m \geq 0$ the polynomial $\sum_{j=0}^{n}(-1)^j q^{mj} \begin{bmatrix} n \\ j \end{bmatrix}_{q^p}$ is divisible by $\prod_{j=1}^{\left\lfloor \frac{n+1}{2} \right\rfloor}\left(1 - q^{\frac{2j-1}{V_p(2j-1)}}\right)$ in $\mathbb{Z}[q]$.

**Conjecture 2.3**

For any odd prime $p$ and any $m \in \mathbb{N}$ there exists a factorization

$$r_n\left(q^{2m+1}, q^{2p}\right) = \sum_{j=0}^{n} q^{(2m+1)j} \begin{bmatrix} n \\ j \end{bmatrix}_{q^{2p}} = \prod_{k=1}^{n}\left(1 + q^{\frac{k}{V_p(k)}}\right) c_p(2m+1, n, q) \quad (2.18)$$

with $c_p(2m+1, n, q) \in \mathbb{Z}[q]$.

**Example**

Let us consider for example $r_n\left(q^5, q^6\right) = \sum_{j=0}^{n} q^{5j} \begin{bmatrix} n \\ j \end{bmatrix}_{q^6}$.

Note that $\left(\frac{k}{V_3(k)}\right)_{k \geq 1} = (1, 2, 1, 4, 5, 2, 7, 8, 1, \cdots)$.

Here we get

$r_0\left(q^5, q^6\right) = 1$, $r_1\left(q^5, q^6\right) = (1+q)c_3(5,1,q)$, $r_2\left(q^5, q^6\right) = (1+q)\left(1+q^2\right)c_3(5,2,q)$,

$r_3\left(q^5, q^6\right) = (1+q)^2\left(1+q^2\right)c_3(5,3,q)$, $r_4\left(q^5, q^6\right) = (1+q)^2\left(1+q^2\right)\left(1+q^4\right)c_3(5,4,q)$,

$r_5\left(q^5, q^6\right) = (1+q)^2\left(1+q^2\right)\left(1+q^4\right)\left(1+q^5\right)c_3(5,5,q), \cdots$.

The first few terms of the sequence $\left(c_3(5, n, q)\right)_{n \geq 0}$ are $1,\ 1-q+q^2-q^3+q^4,\ 1-q+q^4-q^6+q^8,\ \left(1-q+q^2-q^3+q^4\right)\left(1-q+q^4-q^5+q^8-q^9+q^{10}-q^{13}+q^{14}\right),\cdots$

From the definition it is clear that $c_3(5,n,1)=1$. More interesting is the sequence

$$\left(c_3(5,n,-1)\right)_{n\geq 0} = \left(1,\ 5,\ 3,\ 5\cdot 3^2,\ 3^3,\ 5\cdot 3^3,\ 3^4,\ 5\cdot 3^4,\ 3^5,\ 5\cdot 3^7,\cdots\right).$$

Let us more generally compute $c_p(2m+1,n,-1)$. To this end we need a special case of Corollary 2.4 which will be proved later ((2.25) and (2.26)):

$$\lim_{q\to -1} \frac{\sum_{j=0}^{2n} q^{(2m+1)j} \begin{bmatrix} 2n \\ j \end{bmatrix}_{q^{2p}}}{(-q;q)_{2n}} = p^n \tag{2.19}$$

and

$$\lim_{q\to -1} \frac{\sum_{j=0}^{2n+1} q^{(2m+1)j} \begin{bmatrix} 2n+1 \\ j \end{bmatrix}_{q^{2p}}}{(-q;q)_{2n+1}} = p^n(2m+1). \tag{2.20}$$

We have

$$c_p(2m+1,n,-1) = \lim_{q\to -1} \frac{\sum_{j=0}^{n} q^{(2m+1)j} \begin{bmatrix} n \\ j \end{bmatrix}_{q^{2p}}}{\prod_{k=1}^{n}\left(1+q^{\frac{k}{V_p(k)}}\right)} = \lim_{q\to -1} \frac{\sum_{j=0}^{n} q^{(2m+1)j} \begin{bmatrix} n \\ j \end{bmatrix}_{q^{2p}}}{(-q;q)_n} \cdot \frac{(-q;q)_n}{\prod_{k=1}^{n}\left(1+q^{\frac{k}{V_p(k)}}\right)}.$$

It remains to compute $\displaystyle\lim_{q\to -1} \frac{(-q;q)_n}{\prod_{j=1}^{n}\left(1+q^{\frac{j}{V_p(j)}}\right)}$. To this end observe that for $j \not\equiv 0 \bmod p$ we have

$V_p(j)=1$ and therefore $\displaystyle \frac{(-q;q)_n}{\prod_{j=1}^{n}\left(1+q^{\frac{j}{V_p(j)}}\right)} = \prod_{2jp\leq n} \frac{1+q^{2jp}}{1+q^{\frac{2jp}{V_p(2jp)}}} \prod_{(2j+1)p\leq n} \frac{1+q^{(2j+1)p}}{1+q^{\frac{(2j+1)p}{V_p((2j+1)p)}}}.$

For $q\to -1$ the first product converges to 1.

Since $\displaystyle\lim_{q\to -1} \frac{1+q^{(2j+1)p}}{1+q^{\frac{(2j+1)p}{V_p((2j+1)p)}}} = V_p((2j+1)p)$ and $V_p\left(2\cdot 4\cdots 2\left\lfloor\frac{n}{2}\right\rfloor\right) = V_p\left(\left\lfloor\frac{n}{2}\right\rfloor!\right)$

we get $\displaystyle \lim_{q\to -1} \frac{(-q;q)_n}{\prod_{j=1}^{n}\left(1+q^{a(j,p)}\right)} = p^{v_p(n!)-v_p\left(\left\lfloor\frac{n}{2}\right\rfloor!\right)}.$

This gives

$$c_p(1, n, -1) = p^{b(n,p)} \tag{2.21}$$

with

$$b(n, p) = \left\lfloor \frac{n}{2} \right\rfloor + v_p(n!) - v_p\left(\left\lfloor \frac{n}{2} \right\rfloor!\right). \tag{2.22}$$

Therefore we get $c_p(2m+1, 2n, -1) = p^{n+v_p((2n)!)-v_p(n!)}$ and
$c_p(2m+1, 2n+1, -1) = (2m+1)p^{n+v_p((2n+1)!)-v_p(n!)}$.

## 2.3. Limit relations

Let us now consider how fast $\sum_{j=0}^{n}(-1)^j q^{mj} \begin{bmatrix} n \\ j \end{bmatrix}_{q^k}$ converges to $0$ for $q \to 1$.

**Theorem 2.2**

$$\lim_{q \to 1} \frac{\sum_{j=0}^{2n}(-1)^j q^{mj} \begin{bmatrix} 2n \\ j \end{bmatrix}_{q^k}}{(q;q^2)_n} = k^n \tag{2.23}$$

and

$$\lim_{q \to 1} \frac{\sum_{j=0}^{2n+1}(-1)^j q^{mj} \begin{bmatrix} 2n+1 \\ j \end{bmatrix}_{q^k}}{(q;q^2)_{n+1}} = mk^n. \tag{2.24}$$

**Proof**

By (1.20) we have

$$r_{2n}(-q^m, q^k) = \left(1 - (1+q^k)q^{2k(n-1)+m} + q^{2m}\right)r_{2n-2}(-q^m, q^k) - \left(1-q^{k(2n-3)}\right)\left(1-q^{2k(n-1)}\right)q^{2m} r_{2n-4}(-q^m, q^k).$$

Therefore

$$h(n, q) = \frac{r_{2n}(-q^m, q^k)}{(q;q^2)_n}$$

satisfies the recurrence

$$h(n, q) = \frac{1-(1+q^k)q^{2k(n-1)+m}+q^{2m}}{1-q^{2n-1}} h(n-1, q) - \frac{\left(1-q^{k(2n-3)}\right)\left(1-q^{2k(n-1)}\right)q^{2m}}{\left(1-q^{2n-1}\right)\left(1-q^{2n-3}\right)} h(n-2, q)$$

with initial values

$h(0,q) = 1$ and $h(1,q) = [m] - q^{m+k}[m-k]$.

For $q \to 1$ we get

$$h(n,1) = \frac{k(4n-3)}{2n-1} h(n-1,1) - \frac{2k^2(n-1)}{2n-1} h(n-2,1)$$

with initial values $h(0,1) = 1$ and $h(1,1) = k$.

This implies $h(n,1) = k^n$ and thus (2.23).

In the same way $r_{2n+1}(-q^m, q^k)$

satisfies the recurrence

$$r_{2n+1}(-q^m, q^k) = \left(1 - (1+q^k)q^{k(2n-1)+m} + q^{2m}\right) r_{2n-1}(-q^m, q^k)$$
$$- \left(1 - q^{k(2n-2)}\right)\left(1 - q^{k(2n-1)}\right) q^{2m} r_{2n-3}(-q^m, q^k).$$

Therefore

$$H(n,q) = \frac{r_{2n+1}(-q^m, q^k)}{(q;q^2)_{n+1}}$$

satisfies the recurrence

$$H(n,q) = \frac{1 - (1+q^k)q^{k(2n-1)+m} + q^{2m}}{1 - q^{2n+1}} H(n-1,q) - \frac{\left(1 - q^{k(2n-1)}\right)\left(1 - q^{2k(n-1)}\right) q^{2m}}{\left(1 - q^{2n-1}\right)\left(1 - q^{2n+1}\right)} H(n-2,q)$$

with initial values

$$H(0,q) = [m] \quad \text{and} \quad H(1,q) = [m]\frac{[m+k] + q^{2m}[2k-m]}{[3]}.$$

For $q \to 1$ this gives

$$H(n,1) = \frac{k(4n-1)}{2n+1} H(n-1,1) - \frac{2k^2(n-1)}{2n+1} H(n-2,1)$$

with initial values $H(0,1) = m$ and $H(1,1) = km$. This gives (2.24).

The same argument gives

**Corollary 2.4**

$$\lim_{q \to -1} \frac{\sum_{j=0}^{n} q^{(2m+1)j} \begin{bmatrix} 2n \\ j \end{bmatrix}_{q^{2k}}}{(-q;q^2)_n} = (2k)^n \qquad (2.25)$$

and

$$\lim_{q \to -1} \frac{\sum_{j=0}^{n} q^{(2m+1)j} \begin{bmatrix} 2n+1 \\ j \end{bmatrix}_{q^{2k}}}{(-q;q^2)_{n+1}} = (2k)^n (2m+1). \tag{2.26}$$

Computer experiments suggest more facts of this kind. One of my conjectures has been proved by W. Sawin [12].

**Theorem 2.3 (W. Sawin [12])**

*Let $r, m \in \mathbb{N}$. Then*

$$\lim_{q \to 1} \frac{\sum_{j=0}^{2n} (-1)^j q^{rj^2+mj} \begin{bmatrix} 2n \\ j \end{bmatrix}_{q^k}}{(q;q^2)_n} = (k-2r)^n \tag{2.27}$$

**Proof**

Let $f(n,r,m,k) = \sum_{j=0}^{n} (-1)^j q^{rj^2+mj} \begin{bmatrix} n \\ j \end{bmatrix}_{q^k}$.

For $k = 2r$ we get from (1.15)

$$\sum_{j=0}^{n} (-1)^j q^{\frac{k}{2}j^2+mj} \begin{bmatrix} n \\ j \end{bmatrix}_{q^k} = \left(1 - q^{\frac{k}{2}+m}\right)\left(1 - q^{3\frac{k}{2}+m}\right) \cdots \left(1 - q^{(2n-1)\frac{k}{2}+m}\right) = \left(q^{\frac{k}{2}+m}; q^k\right)_n. \tag{2.28}$$

Since $\lim_{q \to 1}(1-q^a) = 0$ and $\lim_{q \to 1} \frac{1-q^a}{1-q^b} = \frac{a}{b}$ we get

$$\lim_{q \to 1} \frac{\sum_{j=0}^{2n} (-1)^j q^{\frac{k}{2}j^2+m} \begin{bmatrix} 2n \\ j \end{bmatrix}_{q^k}}{(q;q^2)_n} = \lim_{q \to 1} \frac{\left(q^{\frac{k}{2}+m}; q^k\right)_{2n}}{(q;q^2)_n} = [n=0] \tag{2.29}$$

and

$$\lim_{q \to 1} \frac{\sum_{j=0}^{2n+1} (-1)^j q^{\frac{k}{2}j^2+mj} \begin{bmatrix} 2n+1 \\ j \end{bmatrix}_{q^k}}{(q;q^2)_{n+1}} = \lim_{q \to 1} \frac{\left(q^{\frac{k}{2}+m}; q^k\right)_{2n+1}}{(q;q^2)_{n+1}} = \left(\frac{k}{2}+m\right)[n=0].$$

Let us consider these limits from another point of view. The denominator of (2.29) has the root $q=1$ with multiplicity $n$. If the limit exists then the numerator also has the root $q=1$ with multiplicity at least $n$. Sawin's trick was to use de L'Hospital's rule combined with Leibniz's rule and the explicit result (2.28).

Observe that

$$\frac{\partial^i f\left(2n, \frac{k}{2}, m, k\right)}{\partial q^i}(1) = \frac{\partial^i \sum_{j=0}^{2n}(-1)^j q^{\frac{k}{2}j^2+mj} \begin{bmatrix} 2n \\ j \end{bmatrix}_{q^k}}{\partial q^i}(1)$$

$$= \sum_{a=0}^{i} \binom{i}{a} \sum_{j=0}^{2n}(-1)^j \frac{\partial^a q^{mj}}{\partial q^a}(1) \frac{\partial^{i-a} q^{\frac{k}{2}j^2} \begin{bmatrix} 2n \\ j \end{bmatrix}_{q^k}}{\partial q^{i-a}}(1)$$

$$= \sum_{a=0}^{i} \binom{i}{a} \sum_{j=0}^{2n}(-1)^j (mj)(mj-1)\cdots(mj-a+1) \frac{\partial^{i-a} q^{\frac{k}{2}j^2} \begin{bmatrix} 2n \\ j \end{bmatrix}_{q^k}}{\partial q^{i-a}}(1)$$

$$= \sum_{a=0}^{i} \binom{i}{a} \sum_{b=0}^{a} [j^b](mj)(mj-1)\cdots(mj-a+1) F(b, i-a, 2n)$$

with

$$F(b, c, 2n) = \sum_{j=0}^{2n}(-1)^j j^b \frac{\partial^c q^{\frac{k}{2}j^2} \begin{bmatrix} 2n \\ j \end{bmatrix}_{q^k}}{\partial q^c}(1).$$

Here we have $b + 2c \leq 2a + 2c \leq 2i$.

It is clear that $F(0, 0, 2n) = \sum_{j=0}^{2n}(-1)^j q^{\frac{k}{2}j^2} \begin{bmatrix} 2n \\ j \end{bmatrix}_{q^k}(1) = [n=0]$.

We want to show that $F(b, c, 2n) = 0$ for $b + c < 2n$.

Let us assume by induction that

$F(b, c, 2n) = 0$ for $b + c < i$ for some fixed $i < 2n$.

Then we get

$$\frac{\partial^i \sum_{j=0}^{2n}(-1)^j q^{\frac{k}{2}j^2+mj}\begin{bmatrix}2n\\j\end{bmatrix}_{q^k}}{\partial q^i}(1) = \sum_{a=0}^{i}\binom{i}{a}[j^a](mj)(mj-1)\cdots(mj-a+1)F(a,i-a,2n) \quad (2.30)$$

$$= \sum_{a=0}^{i}\binom{i}{a}m^a F(a,i-a,2n).$$

Since

$$\sum_{j=0}^{2n}(-1)^j q^{\frac{k}{2}j^2+mj}\begin{bmatrix}2n\\j\end{bmatrix}_{q^k} = \left(q^{\frac{k}{2}+m};q^k\right)_{2n}$$

the polynomial (2.30) vanishes for all $m$ for $i < 2n$. Therefore $F(a, i-a, 2n) = 0$ for each $a$.

Thus $F(b, c, 2n) = 0$ for $b + c = i$.

We thus know that $F(b, c, 2n) = 0$ if $b + c < 2n$. Thus

$$\frac{\partial^n f\left(2n,\frac{k}{2},m,k\right)}{\partial q^n}(1) = \sum_{a=0}^{n}\binom{n}{a}\sum_{j=0}^{2n}(-1)^j (mj)(mj-1)\cdots(mj-a+1)\frac{\partial^{n-a} q^{\frac{k}{2}j^2}\begin{bmatrix}2n\\j\end{bmatrix}_{q^k}}{\partial q^{n-a}}(1)$$

is a linear sum of $F(b, c, 2n)$ with $b + 2c \le 2n$. If $c \ge 1$ then $b + c < 2n$ and therefore $F(b, c, 2n) = 0$. The only remaining case is $F(2n, 0, 2n)$.

In this case we get from the well-known formula $\sum_{k=0}^{n}(-1)^{n-k}k^m\binom{n}{k} = n![n=m]$ that

$$F(2n,0,2n) = \sum_{j=0}^{2n}(-1)^j j^{2n} q^{\frac{k}{2}j^2}\begin{bmatrix}2n\\j\end{bmatrix}_{q^k}(1) = \sum_{j=0}^{2n}(-1)^j j^{2n}\binom{2n}{j} = (2n)!.$$

We shall now use these results in the general case.

We write $f(n, r, m, k)$ in the form

$$f(n,r,m,k) = \sum_{j=0}^{n}(-1)^j q^{rj^2+mj}\begin{bmatrix}n\\j\end{bmatrix}_{q^k} = \sum_{j=0}^{n}(-1)^j q^{\left(r-\frac{k}{2}\right)j^2+mj}\cdot q^{\frac{k}{2}j^2}\begin{bmatrix}n\\j\end{bmatrix}_{q^k}.$$

This gives

$$\frac{\partial^i f(2n,r,m,k)}{\partial q^i}(1) = \sum_{a=0}^{i}\binom{i}{a}\sum_{j=0}^{2n}(-1)^j \frac{\partial^a q^{\left(r-\frac{k}{2}\right)j^2+mj}}{\partial q^a}(1)\frac{\partial^{i-a} q^{\frac{k}{2}j^2}\begin{bmatrix}2n\\j\end{bmatrix}_{q^k}}{\partial q^{i-a}}(1).$$

Now $\dfrac{\partial^{a} q^{\left(r-\frac{k}{2}\right)j^{2}+mj}}{\partial q^{a}}(1)$ is a polynomial of degree $2a$ in $j$. Therefore the right-hand side can be written as a linear combination of the terms $F(b,c,2n)$ where $b+2c \leq 2a+2(i-a)=2i.$

We get as above

$$\dfrac{\partial^{n} f(2n,r,m,k)}{\partial q^{n}}(1) = \sum_{a=0}^{n}\binom{n}{a}\sum_{j=0}^{2n}(-1)^{j}\dfrac{\partial^{a} q^{\left(r-\frac{k}{2}\right)j^{2}+mj}}{\partial q^{a}}(1)\dfrac{\partial^{n-a} q^{\frac{k}{2}j^{2}}\begin{bmatrix}2n \\ j\end{bmatrix}_{q^{k}}}{\partial q^{n-a}}(1)$$

$$= F(2n,0,2n)[j^{2n}]\dfrac{\partial^{n} q^{\left(r-\frac{k}{2}\right)j^{2}+mj}}{\partial q^{n}}(1).$$

Now

$$[j^{2n}]\dfrac{\partial^{n} q^{\left(r-\frac{k}{2}\right)j^{2}+mj}}{\partial q^{n}}(1)$$

$$= [j^{2n}]\left(\left(r-\dfrac{k}{2}\right)j^{2}+mj\right)\left(\left(r-\dfrac{k}{2}\right)j^{2}+mj-1\right)\cdots\left(\left(r-\dfrac{k}{2}\right)j^{2}+mj-n+1\right) = \left(r-\dfrac{k}{2}\right)^{n}.$$

Therefore we finally get

$$\dfrac{\partial^{n} f(2n,r,m,k)}{\partial q^{n}}(1) = (2n)!\left(r-\dfrac{k}{2}\right)^{n}.$$

By Gauss's theorem we have $\left(q;q^{2}\right)_{n} = f(2n,0,0,1).$

Thus

$$\dfrac{\partial^{n} f(2n,0,0,1)}{\partial q^{n}}(1) = (2n)!\left(-\dfrac{1}{2}\right)^{n}.$$

This implies by using de L'Hospital's rule

$$\lim_{q\to 1}\frac{\sum_{j=0}^{2n}(-1)^j q^{rj^2+mj}\begin{bmatrix}2n\\j\end{bmatrix}_{q^k}}{\sum_{j=0}^{2n}(-1)^j\begin{bmatrix}2n\\j\end{bmatrix}_{q^k}}=\lim_{q\to 1}\frac{\dfrac{\partial^n \sum_{j=0}^{2n}(-1)^j q^{rj^2+mj}\begin{bmatrix}2n\\j\end{bmatrix}_{q^k}}{\partial q^n}}{\dfrac{\partial^n \sum_{j=0}^{2n}(-1)^j q^{j}\begin{bmatrix}2n\\j\end{bmatrix}_{q}}{\partial q^n}}=\frac{\left(r-\dfrac{k}{2}\right)^n}{\left(-\dfrac{1}{2}\right)^n}=(k-2r)^n.$$

In an analogous way we get

**Theorem 2.4**

$$\lim_{q\to 1}\frac{\sum_{j=0}^{2n+1}(-1)^j q^{rj^2+mj}\begin{bmatrix}2n+1\\j\end{bmatrix}_{q^k}}{(q;q^2)_{n+1}}=\big((2n+1)r+m\big)(k-2r)^n. \qquad (2.31)$$

**Proof**

Let

$$F(b,c,2n+1)=\sum_{j=0}^{2n+1}(-1)^j j^b \frac{\partial^c q^{\frac{k}{2}j^2}\begin{bmatrix}2n+1\\j\end{bmatrix}_{q^k}}{\partial q^c}(1).$$

As above we get that $F(b,c,2n+1)=0$ if $b+c<2n+1$ and that

$$\frac{\partial^{n+1} f(2n+1,r,m,k)}{\partial q^{n+1}}(1)=\sum_{a=0}^{n+1}\binom{n}{a}\sum_{j=0}^{2n+1}(-1)^j \frac{\partial^a q^{\left(r-\frac{k}{2}\right)j^2+mj}}{\partial q^a}(1)\frac{\partial^{n-a} q^{\frac{k}{2}j^2}\begin{bmatrix}2n+1\\j\end{bmatrix}_{q^k}}{\partial q^{n-a}}(1)$$

is a linear combination of $F(b,c,2n+1)$ with $b+2c\le 2n+2$.

Here for all pairs $(b,c)$ which occur in this sum $F(b,c,2n+1)=0$ except for $(b,c)=(2n+2,0),\ (2n+1,0)$ and $(2n,1)$.

Thus

$$\frac{\partial^{n+1} f(2n+1,r,m,k)}{\partial q^{n+1}}(1) = \sum_{j=0}^{2n+1} (-1)^j \left[ \left[ j^{2n+2} \right] \frac{\partial^{n+1} q^{\left(r-\frac{k}{2}\right)j^2+mj}}{\partial q^{n+1}}(1) \right] j^{2n+2} \binom{2n+1}{j}$$

$$+ \sum_{j=0}^{2n+1} (-1)^j \left[ \left[ j^{2n+1} \right] \frac{\partial^{n+1} q^{\left(r-\frac{k}{2}\right)j^2+mj}}{\partial q^{n+1}}(1) \right] j^{2n+1} \binom{2n+1}{j}$$

$$+ \binom{n+1}{1} \sum_{j=0}^{2n+1} (-1)^j \left[ \left[ j^{2n} \right] \frac{\partial^n q^{\left(r-\frac{k}{2}\right)j^2+mj}}{\partial q^n}(1) \right] j^{2n} \frac{\partial \left( q^{\frac{kj^2}{2}} \left[ \begin{matrix} 2n+1 \\ j \end{matrix} \right]_{q^k} \right)}{\partial q}(1).$$

For the first term we get

$$\left[ j^{2n+2} \right] \frac{\partial^{n+1} q^{\left(r-\frac{k}{2}\right)j^2+mj}}{\partial q^{n+1}}(1) = \left( r - \frac{k}{2} \right)^{n+1} \quad \text{and} \quad \sum_{j=0}^{2n+1} (-1)^j j^{2n+2} \binom{2n+1}{j} = -(n+1)(2n+1)(2n+1)!.$$

For the second one we have

$$\left[ j^{2n+1} \right] \frac{\partial^{n+1} q^{\left(r-\frac{k}{2}\right)j^2+mj}}{\partial q^{n+1}}(1) = \left[ j^{2n+1} \right] \left( \left( r - \frac{k}{2} \right) j^2 + mj \right) \left( \left( r - \frac{k}{2} \right) j^2 + mj - 1 \right) \cdots \left( \left( r - \frac{k}{2} \right) j^2 + mj - n \right)$$

$$= \left( r - \frac{k}{2} \right)^n (n+1)m$$

and

$$\sum_{j=0}^{2n+1} (-1)^j j^{2n+1} \binom{2n+1}{j} = -(2n+1)!.$$

In the last term we have $\left[ j^{2n} \right] \dfrac{\partial^n q^{\left(r-\frac{k}{2}\right)j^2+mj}}{\partial q^n}(1) = \left( r - \dfrac{k}{2} \right)^n$.

There remains $\displaystyle\sum_{j=0}^{2n+1} (-1)^j j^{2n} \frac{\partial \left( q^{\frac{kj^2}{2}} \left[ \begin{matrix} 2n+1 \\ j \end{matrix} \right]_{q^k} \right)}{\partial q}(1).$

Note that $\dfrac{\partial \left[ \begin{matrix} n \\ j \end{matrix} \right]_{q^k}}{\partial q}(1) = \binom{j+1}{2} \binom{n}{j+1} k.$ This is easily verified by induction. It is obviously true for $j = 0$ and all $n$. Let us suppose that it holds for $j-1$ and all $n$. It is obviously true for $j$ and $n \leq j$. Let it hold for $n$. Then we get

$$\frac{\partial \begin{bmatrix} n+1 \\ j \end{bmatrix}_{q^k}}{\partial q}(1) = \frac{\partial \left( q^{jk} \begin{bmatrix} n \\ j \end{bmatrix}_{q^k} + \begin{bmatrix} n \\ j-1 \end{bmatrix}_{q^k} \right)}{\partial q} = jk\binom{n}{j} + \binom{j+1}{2}\binom{n}{j+1}k + \binom{j}{2}\binom{n}{j}k$$

$$= \binom{j+1}{2}\left(\binom{n}{j+1} + \binom{n}{j}\right)k = \binom{j+1}{2}\binom{n+1}{j+1}k.$$

Therefore we get

$$\frac{\partial \left( q^{\frac{kj^2}{2}} \begin{bmatrix} 2n+1 \\ j \end{bmatrix}_{q^k} \right)}{\partial q}(1) = \binom{j+1}{2}\binom{2n+1}{j+1}k + \frac{j^2 k}{2}\binom{2n+1}{j} = \frac{(2n+1)jk}{2}\binom{2n+1}{j}.$$

This gives

$$\sum_{j=0}^{2n+1}(-1)^j j^{2n} \frac{\partial \left( q^{\frac{kj^2}{2}} \begin{bmatrix} 2n+1 \\ j \end{bmatrix}_{q^k} \right)}{\partial q}(1) = \frac{(2n+1)k}{2}\sum_{j=0}^{2n+1}(-1)^j j^{2n+1}\binom{2n+1}{j} = \frac{(2n+1)k}{2}(2n+1)!.$$

Combining all results we get

$$\frac{\partial^{n+1} f(2n+1, r, m, k)}{\partial q^{n+1}}(1) = \left(r - \frac{k}{2}\right)^n (2n+1)!\left(-(n+1)(2n+1)\left(r - \frac{k}{2}\right) - (n+1)m + (n+1)(2n+1)k/2\right)$$

$$= -(n+1)\left(r - \frac{k}{2}\right)^n (2n+1)!\left((2n+1)r + m\right).$$

Now observe that $(q; q^2)_{n+1} = f(2n+1, 0, 1, 1)$.

This implies

$$\frac{\partial^{n+1} f(2n+1, 0, 1, 1)}{\partial q^{n+1}}(1) = -(n+1)\left(-\frac{1}{2}\right)^n (2n+1)!.$$

Therefore

$$\lim_{q \to 1} \frac{\sum_{j=0}^{2n+1}(-1)^j q^{rj^2 + mj} \begin{bmatrix} 2n+1 \\ j \end{bmatrix}_{q^k}}{(q; q^2)_{n+1}} = \lim_{q \to 1} \frac{\dfrac{\partial^{n+1} \sum_{j=0}^{2n+1}(-1)^j q^{rj^2 + mj} \begin{bmatrix} 2n+1 \\ j \end{bmatrix}_{q^k}}{\partial q^{n+1}}}{\dfrac{\partial^{n+1} \sum_{j=0}^{2n+1}(-1)^j q^j \begin{bmatrix} 2n+1 \\ j \end{bmatrix}_q}{\partial q^{n+1}}}$$

$$= \frac{-(n+1)\left(r-\frac{k}{2}\right)^{n}(2n+1)!((2n+1)r+m)}{-(n+1)\left(-\frac{1}{2}\right)^{n}(2n+1)!} = (k-2r)^{n}((2n+1)r+m).$$

Analogous results hold for binomials of the form $\sum_{j=0}^{n} q^{rj^{2}+mj} \begin{bmatrix} n \\ j \end{bmatrix}_{q^{k}}$. It is easy to check that

$$\lim_{q\to -1} \sum_{j=0}^{n} q^{rj^{2}+mj} \begin{bmatrix} n \\ j \end{bmatrix}_{q^{k}} = 0 \text{ for } n>0 \text{ if and only if } r+m \equiv 1 \bmod 2 \text{ and } k \equiv 0 \bmod 2.$$

Therefore we get

**Corollary 2.5**

*If $r+m \equiv 1 \bmod 2$ then*

$$\lim_{q\to -1} \frac{\sum_{j=0}^{2n} q^{rj^{2}+mj} \begin{bmatrix} 2n \\ j \end{bmatrix}_{q^{2k}}}{(-q;q)_{2n}} = (k-r)^{n} \tag{2.32}$$

*and*

$$\lim_{q\to -1} \frac{\sum_{j=0}^{2n+1} q^{rj^{2}+mj} \begin{bmatrix} 2n+1 \\ j \end{bmatrix}_{q^{2k}}}{(-q;q)_{2n+1}} = ((2n+1)r+m)(k-r)^{n}. \tag{2.33}$$

This follows from (2.27) because $(-q)^{rj^{2}+mj} = (-1)^{j} q^{rj^{2}+mj}$ and $\lim_{q\to -1} \frac{(-q;q)_{2n}}{(-q;q^{2})_{n}} = 2^{n}$.

### 3. Some Hankel determinants

**3.1.** There is a close connection with Hankel determinants.

**Proposition 3.1**

*For each $m \in \mathbb{N}$ the sequences $\left(f\left(n, q^{2m+1}, q\right)\right)_{n\geq 0}$, $\left(F\left(2n, q^{m}, q\right)\right)_{n\geq 0}$, $\left(F\left(2n+1, q^{m}, q\right)\right)_{n\geq 0}$ and $\left(F\left(n, q^{m}, q\right)\right)_{n\geq 0}$ satisfy homogeneous linear recurrences with constant coefficients.*

**Proof**

By (1.15) we see that $\sum_{j=0}^{r}(-1)^{j}q^{\binom{j}{2}}\begin{bmatrix}r\\j\end{bmatrix}q^{m(n-j)}=q^{mn}\sum_{j=0}^{r}(-1)^{j}q^{\binom{j}{2}}\begin{bmatrix}r\\j\end{bmatrix}q^{-mj}=q^{mn}\left(\frac{1}{q^{m}};q\right)_{r}=0$

for $r > m$.

By Corollary 2.1 $f(n,q^{2m+1},q)$ is a linear sum of terms $q^{jn}$ with $0 \leq j < m+1$ and therefore satisfies the recurrence

$$\sum_{j=0}^{m+1}(-1)^{j}q^{\binom{j}{2}}\begin{bmatrix}m+1\\j\end{bmatrix}f(n-j,q^{2m+1},q)=0.$$

The same argument applies to the other cases.

If a sequence $(x(n))_{n\geq 0}$ satisfies a linear homogeneous recurrence of order $r$ with constant coefficients then the rows and columns of the Hankel matrix $(x(i+j))_{i,j=0}^{n}$ are linear dependent for $n \geq r$ and therefore their determinants vanish.

Thus we get

**Corollary 3.1**

*For each $m \in \mathbb{Z}$ the Hankel determinants $d(n,q^{2m+1},q) = \det\left(f(i+j,q^{2m+1},q)\right)_{i,j=0}^{n}$,*

$D_{0}(n,q^{m},q) = \det\left(F(2i+2j,q^{m},q)\right)_{i,j=0}^{n}$, $D_{1}(n,q^{m},q) = \det\left(F(2i+2j+1,q^{m},q)\right)_{i,j=0}^{n}$ *and*

$D(n,q^{m},q) = \det\left(F(i+j,q^{m},q)\right)_{i,j=0}^{n}$ *vanish for all sufficiently large numbers n.*

### 3.2. Orthogonal polynomials and Hankel determinants

For readers who are not familiar with orthogonal polynomials and Hankel determinants I first recall some well-known facts (cf. e.g. [4]).

A sequence of real-valued monic polynomials $(p_n(x))_{n\geq 0}$ with $\deg p_n = n$ is called orthogonal with respect to a linear functional $F$ if $F(p_n(x)p_m(x)) = 0$ for $m \neq n$ and $F(p_n(x)^2) \neq 0$.

In particular for $m = 0$ we get

$$F(p_n(x)) = [n=0]. \tag{2.34}$$

Since $\deg p_n = n$ the linear functional $F$ is uniquely determined by (2.34).

On the other hand $F$ is also uniquely determined by the moments $a(n) = F(x^n)$ for $n \in \mathbb{N}$.

By Favard's theorem there exist numbers $\sigma(n), \tau(n)$ such that

$$p_n(x) = (x - \sigma(n-1))p_{n-1}(x) - \tau(n-2)p_{n-2}(x). \tag{2.35}$$

On the other hand for given sequences $\sigma(n)$ and $\tau(n) \neq 0$ formula (2.35) together with the initial values $p_0(x) = 1$ and $p_{-1}(x) = 0$ defines a sequence of polynomials which are orthogonal with respect to the linear functional $F$ defined by (2.34). If $F(x^n) = a(n)$ then

$$p_n(x) = \frac{1}{\det(a(i+j))_{i,j=0}^{n-1}} \det \begin{pmatrix} a(0) & a(1) & \cdots & a(n-1) & 1 \\ a(1) & a(2) & \cdots & a(n) & x \\ a(2) & a(3) & \cdots & a(n+1) & x^2 \\ \vdots & & & & \vdots \\ a(n) & a(n+1) & \cdots & a(2n-1) & x^n \end{pmatrix}. \tag{2.36}$$

If we define $a(n, j)$ by

$$\begin{aligned} a(0, j) &= [j = 0] \\ a(n, 0) &= \sigma(0)a(n-1, 0) + \tau(0)a(n-1, 1) \\ a(n, j) &= a(n-1, j-1) + \sigma(j)a(n-1, j) + \tau(j)a(n-1, j+1), \end{aligned} \tag{2.37}$$

then $a(n, 0) = a(n)$ and the Hankel determinant $\det(a(i+j))_{i,j=0}^n$ is given by

$$\det(a(i+j))_{i,j=0}^n = \prod_{i=1}^n \prod_{j=0}^{i-1} \tau(j). \tag{2.38}$$

Thus if all $\tau(j)$ have a closed expression so does the Hankel determinant.

Note that (2.35) immediately implies that the polynomials $\overline{p}_n(x) = r^n p_n\left(\frac{x}{r}\right)$ are also orthogonal and satisfy

$$\overline{p}_n(x) = (x - r\sigma(n-1))\overline{p}_{n-1}(x) - r^2\tau(n-2)\overline{p}_{n-2}(x). \tag{2.39}$$

Their moments are given by

$$F((rx)^n) = r^n a(n). \tag{2.40}$$

This scenario can be used to guess and in some cases also to prove the Hankel determinants of a given sequence $(a(n))$.

Let me demonstrate this method in the known case (cf. [6]) of the sequence of Rogers-Szegö polynomials $\left(r_n(s,q)\right)_{n\geq 0} = \left(\sum_{j=0}^{n}\begin{bmatrix}n\\j\end{bmatrix}_q s^j\right)_{n\geq 0}$. Using (2.36) we compute the polynomials $p_n(x,s,q)$ for small values of $n$ and determine $\sigma(n,s,q)$ and $\tau(n,s,q)$ using (2.35). Then we guess that

$$\sigma(n,s,q) = q^n(1+s) \tag{2.41}$$

and

$$\tau(n,s,q) = q^n s\left(q^{n+1}-1\right). \tag{2.42}$$

By (2.38) the Hankel determinants become

$$d(n,s,q) = \det\left(r_{i+j}(s,q)\right)_{i,j=0}^{n} = q^{\binom{n+1}{3}}(-s)^{\binom{n+1}{2}}\prod_{j=1}^{n}(q;q)_j. \tag{2.43}$$

Till now (2.43) is only a guess. In order to prove it we compute $a(n,j) = a(n,j,s,q)$ by (2.37).

Again we guess that

$$a(n,j,s,q) = r_{n-j}(s,q)\begin{bmatrix}n\\j\end{bmatrix}_q. \tag{2.44}$$

Now it is easy to verify (2.37) and we are done.

The same procedure can be applied to obtain

$$p_n(x,s,q) = \sum_{j=0}^{n}(-1)^j q^{\binom{j}{2}}\begin{bmatrix}n\\j\end{bmatrix} r_j\left(s,q^{-1}\right)x^{n-j} \tag{2.45}$$

by using Favard's theorem (2.35).

**Remark**

Comparing with [8], 14.24, we see that the polynomials $p_n(x,s,q)$ are the

Al-Salam-Carlitz I polynomials $p_n(x,s,q) = U_n^{(s)}(x;q) = \sum_{k=0}^{n}\begin{bmatrix}n\\k\end{bmatrix}(-1)^k q^{\binom{k}{2}}s^k\prod_{j=0}^{n-k-1}(x-q^j)$.

The above result is therefore equivalent with the fact that the Rogers-Szegö polynomials $r_n(s,q) = \sum_{j=0}^{n}\begin{bmatrix}n\\j\end{bmatrix}_q s^j$ are the moments of the Al-Salam-Carlitz I polynomials with respect to the linear functional $F$. This is a $q$-analogue of the trivial fact that the linear functional $F$ defined by $F(f) = f(s+1)$ satisfies $F\left((x-1-s)^n\right) = [n=0]$ and $F(x^n) = (s+1)^n$.

Since no Hankel determinant vanishes there is no $s \neq 0$ for which the Rogers-Szegö polynomials satisfy a recursion with constant coefficients.

**3.3. Hankel determinants of** $f(n,s,q) = \dfrac{\sum_{j=0}^{n} s^j \begin{bmatrix} n \\ j \end{bmatrix}_{q^2}}{(-q;q)_n}$.

By (2.3) $f(n,s,q) = \sum_{j=0}^{n} (-1)^j \begin{bmatrix} n \\ j \end{bmatrix} \dfrac{\prod_{i=0}^{j-1}(q^{2i+1} - s)}{(-q;q)_j}$.

Let us consider more generally

$$h(n,s,t,q) = \sum_{j=0}^{n} (-1)^j \begin{bmatrix} n \\ j \end{bmatrix}_q \dfrac{\prod_{i=0}^{j-1}(q^{2i+1}t - s)}{(-qt;q)_j}. \tag{2.46}$$

Note that besides $h(n,s,1,q) = f(n,s,q)$ we also have $h(n,s,0,q) = r_n(s,q) = \sum_{j=0}^{n} \begin{bmatrix} n \\ j \end{bmatrix}_q s^j$

and $h(n,0,t,q) = \dfrac{1}{(-qt;q)_n}$. The last identity is known as Cauchy's identity.

The sequence $(h(n,s,t,q))$ satisfies the recurrence

$$h(n,s,t,q) = \dfrac{(1+s)h(n-1,s,t,q) + (q^{n-1} - 1)sh(n-2,s,t,q)}{1+q^n t}$$

with initial values $h(0,s,t,q) = 1$ and $h(1,s,t,q) = \dfrac{1+s}{1+qt}$.

In the same way as in Lemma 1.3 by applying qZeil we verify the identity

$$h(n,s,t,q) = \dfrac{\sum_{j=0}^{\lfloor \frac{n}{2} \rfloor} (-1)^j q^{j^2} (q;q^2)_j \begin{bmatrix} n \\ 2j \end{bmatrix} (st)^j r_{n-2j}(s,q)}{(-qt;q)_n}.$$

For $t = 1$ this reduces to (1.35).

This also leads to the generating function

$$\sum_{n \geq 0} \dfrac{h(n,s,t,q)}{(q;q)_n} z^n = \dfrac{e_q(sz)e_q(z)}{e_{q^2}(qstz^2)}. \tag{2.47}$$

It is easy to guess that

$$\sigma(n,s,t) = q^n(1+s)\frac{1+q^{n-1}(1+q)t - q^{2n}t}{(1+q^{2n-1}t)(1+q^{2n+1}t)} \tag{2.48}$$

and

$$\tau(n,s,t) = -\frac{q^n(1+q^nt)(1-q^{n+1})}{(1+q^{2n}t)(1+q^{2n+1}t)^2(1+q^{2n+2}t)}(s-q^{2n+1}t)(1-q^{2n+1}st). \tag{2.49}$$

This implies

$$d(n,s,q) = q^{\frac{n^2(n+1)}{2}}\frac{\prod_{j=0}^{n}\left(\frac{s}{q^{2j-1}};q^2\right)_{2j}\prod_{j=0}^{n}(q;q)_j}{\prod_{j=0}^{n}(-q;q)_{j+n}}. \tag{2.50}$$

Finally we also guess that

$$p(n,x,s,t,q) = \sum_{j=0}^{n}(-1)^j q^{\binom{j}{2}}\begin{bmatrix}n\\j\end{bmatrix}h\left(j,s,q^{2n}t,\frac{1}{q}\right)x^{n-j}$$

$$= \sum_{j=0}^{n}q^{\binom{n-j}{2}}\begin{bmatrix}n\\j\end{bmatrix}\frac{\prod_{i=0}^{n-j-1}(q^{2j+2i+1}t-s)}{(-q^{n+j}t;q)_{n-j}}\prod_{i=0}^{j-1}(x-q^i) \tag{2.51}$$

and

$$a(n,k,s,t) = \begin{bmatrix}n\\k\end{bmatrix}h(n-k,s,q^{2k}t,q). \tag{2.52}$$

To give a formal proof of (2.50) it would suffice to verify (2.37) for these polynomials.

As Dennis Stanton kindly informed me these results already follow from his paper [7] with Mourad E. H. Ismail. Instead of Rogers-Szegö polynomials they consider the intimately related continuous $q$–Hermite polynomials $H_n(x|q)$. These satisfy $H_{-1}(x|q) = 0$, $H_0(x|q) = 1$ and

$$2xH_n(x|q) = H_{n+1}(x|q) + (1-q^n)H_{n-1}(x|q) \tag{2.53}$$

for $n > 0$.

The connection is given by the well-known formula

$$H_n(x|q) = \sum_{k=0}^{n} \begin{bmatrix} n \\ k \end{bmatrix} e^{i(n-2k)\theta} = e^{in\theta} r_n(e^{-2i\theta}, q) \quad \text{for } x = \cos\theta, \tag{2.54}$$

which can be verified by comparing the recursions of both sides.

In order to state their results we need the monic version

$$p_n(x) = p_n(x; a, b, c; q) = \frac{(qa; q)_n (qc; q)_n}{(q^{n+1} ab; q)_n} P_n(x; a, b, c; q) \tag{2.55}$$

of the big $q$ – Jacobi polynomials $P_n(x; a, b, c; q)$. They satisfy ( cf. [8], (14.5.4))

$$x p_n(x) = p_{n+1}(x) + (1 - (A_n + C_n)) p_n(x) + A_{n-1} C_n p_{n-1}(x) \tag{2.56}$$

with

$$A_n = \frac{(1 - aq^{n+1})(1 - abq^{n+1})(1 - cq^{n+1})}{(1 - abq^{2n+1})(1 - abq^{2n+2})}$$

$$C_n = -acq^{n+1} \frac{(1 - q^n)(1 - abc^{-1} q^n)(1 - bq^n)}{(1 - abq^{2n})(1 - abq^{2n+1})}. \tag{2.57}$$

In [7], Theorem 4.1. (D) Ismail and Stanton give a representation of the polynomials $\dfrac{H_n(x|q^2)}{(-q;q)_n}$ as moments of some explicitly given $q$ – integral. In our terminology their result can be stated in the following way:

Let $F$ be the linear functional which satisfies $F(p_n) = [n = 0]$ for the monic big $q$ – Jacobi polynomials $p_n\left(x; -\dfrac{1}{s\sqrt{q}}, \dfrac{s}{\sqrt{q}}, -\dfrac{s}{\sqrt{q}}; q\right)$.

In this case we get (cf. [8], 14.5)

$$A_n C_{n+1} = \frac{q^{n+1}(1 + q^n)(1 - q^{n+1})(1 - q^{2n+1} s^2)(q^{2n+1} - s^2)}{s^2 (1 + q^{2n})(1 + q^{2n+1})^2 (1 + q^{2n+2})}$$

and

$$A_n + C_n - 1 = -\frac{1 + s^2}{s} q^{n + \frac{1}{2}} \frac{1 + q^{n-1}(1 + q) - q^{2n}}{(1 + q^{2n-1})(1 + q^{2n+1})}.$$

[7], Theorem 4.1.(D) implies

$$F(x^n) = \left(-\sqrt{q}\right)^n \frac{H_n\left(\frac{1}{2}\left(s+\frac{1}{s}\right)\mid q^2\right)}{(-q;q)_n}. \tag{2.58}$$

By (2.39) and

$$H_n\left(\frac{1}{2}\left(s+\frac{1}{s}\right)\mid q\right) = s^{-n} r_n(s^2, q) \tag{2.59}$$

this means that the monic orthogonal polynomials corresponding to $\dfrac{r_n(s^2, q^2)}{(-q;q)_n}$ are

$$\left(-\frac{s}{\sqrt{q}}\right)^n p_n\left(-\frac{\sqrt{q}}{s}x; -\frac{1}{s\sqrt{q}}, \frac{s}{\sqrt{q}}, -\frac{s}{\sqrt{q}}; q\right)$$

in which case we get $A_n + C_n - 1 = q^n (1+s^2) \dfrac{1 + q^{n-1}(1+q) - q^{2n}}{(1+q^{2n-1})(1+q^{2n+1})}$ and

$$A_n C_{n+1} = \frac{q^n (1+q^n)(1-q^{n+1})(1-q^{2n+1}s^2)(q^{2n+1} - s^2)}{(1+q^{2n})(1+q^{2n+1})^2(1+q^{2n+2})}.$$

If we now replace $s^2 \to s$ we get (2.48) and (2.49) for $t = 1$.

Thus we get

**Theorem 3.1 (M.E.H. Ismail and D. Stanton [7])**

*The Hankel determinants $d(n, s, q)$ of $f(n, s, q) = \dfrac{\sum_{j=0}^n s^j \begin{bmatrix} n \\ j \end{bmatrix}_{q^2}}{(-q;q)_n}$ are given by*

$$d(n, s, q) = q^{\frac{n^2(n+1)}{2}} \frac{\prod_{j=0}^n \left(\dfrac{s}{q^{2j-1}}; q^2\right)_{2j} \prod_{j=0}^n (q;q)_j}{\prod_{j=0}^n (-q;q)_{j+n}}. \tag{2.60}$$

## 3.4. Hankel determinants of $F(2n,s,q)$ and $F(2n+1,s,q)$.

Let us consider

$$H(n,s,t,q) = \sum_{j=0}^{n} \begin{bmatrix} n \\ j \end{bmatrix}_{q^2} \frac{\prod_{i=0}^{2j-1}(q^i t - s)}{(qt^2;q^2)_j} \tag{2.61}$$

which satisfies $H(n,s,1,q) = F(2n,s,q)$, $H(n,s,q,q) = \dfrac{1-q}{1-s} F(2n+1,s,q)$ and

$H(n,s,0,q) = \sum_{j=0}^{n} s^{2j} \begin{bmatrix} n \\ j \end{bmatrix}_{q^2}$. We have also $H(n,0,t,q) = \dfrac{1}{(qt^2;q^2)_n}$.

The sequence $(H(n,s,t,q))$ satisfies the recurrence

$$H(n,s,t,q) = \frac{\left(1+s^2 - q^{2n-2}(1+q)st\right)H(n-1,s,t,q) + \left(1-q^{2n-2}\right)s^2 H(n-2,s,t,q)}{1-q^{2n-1}t^2}$$

with initial values $H(0,s,t,q) = 1$ and $H(1,s,t,q) = \dfrac{1+(1+q)s+s^2}{1-qt^2}$.

We also get

$$H(n,s,t,q) = \frac{\sum_{j=0}^{n}(-1)^j q^{\binom{j}{2}} \begin{bmatrix} n \\ j \end{bmatrix}_{q^2} s^j t^j r_{n-j}(s^2,q^2)}{(qt^2;q^2)_n}, \tag{2.62}$$

which for $t=1$ reduces to (1.29).

For the generating function we get

$$\sum_{n\geq 0} \frac{H(n,s,t,q)}{(q^2;q^2)_n} z^n = \frac{e_{q^2}(s^2 z) e_{q^2}(z)}{e_q(stz)}. \tag{2.63}$$

Here we get

$$\sigma(n,s,t) = q^{2n-2} \frac{(1+s^2)q^2\left(1 - q^{2n-1}t^2 - q^{2n-3}t^2 + q^{4n-1}t^2\right) + st(1+q)\left(1 - q^{2n} - q^{2n+2} + q^{4n-1}t^2\right)}{\left(1-q^{4n+1}t^2\right)\left(1-q^{4n-3}t^2\right)}$$
(2.64)

and

$$\tau(n,s,t) = \tau(n,0,t)\left(s - q^{2n}t\right)\left(s - q^{2n+1}t\right)\left(1 - q^{2n}st\right)\left(1 - q^{2n+1}st\right) \tag{2.65}$$

with

$$\tau(n,0,t) = -q^{2n} \frac{(1-q^{2n-1}t^2)(1-q^{2n+2})}{(1-q^{4n-1}t^2)(1-q^{4n+1}t^2)^2(1-q^{4n+3}t^2)}. \tag{2.66}$$

In this case we get

$$p_n(x,s,t) = \sum_{j=0}^{n} (-1)^j q^{2\binom{j}{2}} \begin{bmatrix} n \\ j \end{bmatrix}_{q^2} H\left(j,s,q^{2n-1}t,\frac{1}{q}\right)$$

$$= \sum_{j=0}^{n} \begin{bmatrix} n \\ j \end{bmatrix}_{q^2} q^{2\binom{n-j}{2}} \frac{\prod_{i=0}^{2n-2j-1}(s-q^{2j+i}t)}{(q^{2n-1+2j}t^2;q^2)_{n-j}} \prod_{i=0}^{j-1}(x-q^{2i}). \tag{2.67}$$

and

$$a(n,k,s,t) = \begin{bmatrix} n \\ k \end{bmatrix}_{q^2} H\left(n-k,s,q^{2k}t,q\right). \tag{2.68}$$

A proof can again be deduced from [7]. At the end of section 4 two formulae are stated, the first one of which is equivalent with

$$F(x^n) = (-q)^n \frac{H_{2n}\left(\frac{1}{2}\left(s+\frac{1}{s}\right)|q\right)}{(q;q^2)_n} \tag{2.69}$$

where $F$ is the linear functional defined by $F(p_n) = [n=0]$ for

$$p_n = p_n\left(x; -\frac{1}{qs^2}, -\frac{s^2}{q^2}, -\frac{s^2}{q}; q^2\right).$$

In this case we get

$$1-(A_n+C_n) = q^{2n-1} \frac{(1+s^4)q^2\left(1-q^{2n-1}-q^{2n-3}+q^{4n-1}\right) - s^2(1+q)\left(1-q^{2n}-q^{2n+2}+q^{4n-1}\right)}{s^2\left(1-q^{4n+1}\right)\left(1-q^{4n-3}\right)}$$

and

$$A_n C_{n+1} = -q^{2n+2} \frac{(1-q^{2n-1})(1-q^{2n+2})}{(1-q^{4n-1})(1-q^{4n+1})^2(1-q^{4n+3})} \frac{\left(s^2+q^{2n}\right)\left(s^2+q^{2n+1}t\right)\left(1+q^{2n}s^2\right)\left(1+q^{2n+1}s^2\right)}{s^4}.$$

The corresponding terms for $F(2n,s,q)$ are $\left(-\frac{s^2}{q}(A_n+C_n-1)\right)_{s^2 \to -s}$ which is (2.64) for

$t=1$ and $\left(\frac{s^4}{q^2}A_n C_{n+1}\right)_{s^2 \to -s}$ which gives (2.65) for $t=1$.

The second formula at the end of [7], Section 4, is equivalent with

$$F(x^n) = (-q^2)^n \frac{s}{1+s^2} \frac{H_{2n+1}\left(\frac{1}{2}\left(s+\frac{1}{s}\right)\mid q\right)}{(q^3;q^2)_n}$$

where $F$ is the linear functional defined by $F(p_n) = [n=0]$ for

$$p_n = p_n\left(x; -\frac{1}{s^2}, -\frac{s^2}{q}, -s^2; q^2\right).$$

In the same way as above we get $\sigma(n,s,q)$ and $\tau(n,s,q)$ of formulae (2.64) - (2.66).

Thus we get

**Theorem 3.2 (M.E.H. Ismail and D. Stanton [7])**

*The Hankel determinants* $D_0(n,s,q) = \det\left(F(2i+2j,s,q)\right)_{i,j=0}^n$ *and*

$D_1(n,s,q) = \det\left(F(2i+2j+1,s,q)\right)_{i,j=0}^n$ *satisfy*

$$D_0(n,s,q) = D_0(n,0,q)\prod_{j=0}^n \left((s-q^{2j})(s-q^{2j+1})\left(s-\frac{1}{q^{2j}}\right)\left(s-\frac{1}{q^{2j+1}}\right)\right)^{n-j}$$

$$D_1(n,s,q) = D_1(n,0,q)(1-s)^{n+1}\prod_{j=0}^n \left((s-q^{2j+2})(s-q^{2j+1})\left(s-\frac{1}{q^{2j+2}}\right)\left(s-\frac{1}{q^{2j+1}}\right)\right)^{n-j}.$$

(2.70)

### 3.5. Hankel determinants of $F(n,s,q)$.

Finally we want to determine the Hankel determinants $D(n,s,q) = \det\left(F(i+j,s,q)\right)_{i,j=0}^n$.

Here we get

$$\sigma(n,s) = (-1)^n \frac{(1-s)^2 u(n) + sv(n)}{(1-s)(1-q^{2n-1})(1-q^{2n+1})} \tag{2.71}$$

with

$$\begin{aligned} u(2n) &= q^{2n}\left(1-q^{2n-1}+q^{2n+1}-q^{4n+1}\right), \\ u(2n+1) &= q^{2n+2}\left(1-q^{4n+1}+q^{2n-1}-q^{2n+1}\right) \end{aligned} \tag{2.72}$$

and

$$v(2n) = \left(1-q^{4n-1}\right)\left(1-q^{2n+1}\right)\left(1-q^{2n}\right)$$
$$v(2n+1) = \left(1-q^{4n+3}\right)\left(1-q^{2n+1}\right)\left(1-q^{2n}\right).$$

The corresponding $\tau(n,s)$ are

$$\tau(2n,s) = \frac{-q^{4n+1}\left(1-\dfrac{s}{q^{2n}}\right)\left(1-\dfrac{s}{q^{2n+1}}\right)\left(1-q^{2n}s\right)\left(1-q^{2n+1}s\right)}{(1-s)^2\left(1-q^{4n+1}\right)^2}, \tag{2.73}$$

$$\tau(2n+1,s) = (1-s)^2 q^{2n+2} \frac{\left(1-q^{2n+1}\right)\left(1-q^{2n+2}\right)}{\left(1-q^{4n+3}\right)^2}. \tag{2.74}$$

This conjecture has been proved by D. Stanton [14]:

**Theorem 3.3**

*The Hankel determinants $D(n,s,q) = \det\left(F(i+j,s,q)\right)_{i,j=0}^{n}$ satisfy*

$$D(n,s,q) = D(n,0,q) \frac{\prod_{j=0}^{\lfloor\frac{n}{2}\rfloor}\left(\left(1-q^{2j}s\right)\left(1-q^{2j+1}s\right)\left(1-\dfrac{s}{q^{2j}}\right)\left(1-\dfrac{s}{q^{2j+1}}\right)\right)^{n-2j}}{(1-s)^{2\lfloor\frac{n+1}{2}\rfloor}} \tag{2.75}$$

*with*

$$D(n,0,q) = (-1)^{\lfloor\frac{(n+1)^2}{4}\rfloor} q^{n^2+n\sum_{j=1}^{n-1}\lfloor\frac{j}{2}\rfloor} \prod_{j=0}^{\lfloor\frac{n-1}{2}\rfloor}\left(\frac{1}{\left(1-q^{4j+1}\right)}\right)^{2n-2j}\left(\frac{\left(1-q^{2j+1}\right)\left(1-q^{2j+2}\right)}{\left(1-q^{4j+3}\right)^2}\right)^{n-1-2j}. \tag{2.76}$$

## 4. Conclusion

B. A. Kupershmidt [9] pointed out the fact that the *alternating* Gauss-like sums $r_n\left(-q^m,q\right)$ can be effectively calculated for *integers* $m \in \mathbb{Z}$ whereas the *non-alternating* sums $r_n\left(q^r,q\right)$ can be effectively calculated for *half-integers* $r \in \dfrac{1}{2}+\mathbb{Z}$. The facts about the Hankel determinants show that there are no other special values $s$ such that $r_n(s,q)$ has similar properties. But there are interesting $q$ – analogues of $(1\pm1)^n$ of the form $r_n\left(\pm q^m, q^k\right)$, which have a factorization $r_n\left(\pm q^m, q^k\right) = b_n(q)u_n(q)$ into a "nice part" $b_n(q)$ (which is a product of terms $1\pm q^\ell$) and an "ugly part", which has nice values for $q=\pm 1$. It seems that we have only scratched the tip of an iceberg and there remains much to be done.